\documentclass[12pt]{article}
\usepackage[]{lineno}

\usepackage{mathrsfs}
\usepackage{amsmath,extpfeil}
\usepackage{stmaryrd}
\usepackage{mathrsfs}
\usepackage{url}
\usepackage{indentfirst}
\usepackage{enumerate}
\usepackage{amsmath,amsfonts,amssymb,amsthm}
\usepackage{amsmath,amssymb,amsthm,amscd}
\usepackage{graphicx,mathrsfs}
\usepackage{appendix}
\usepackage[numbers,sort&compress]{natbib}
\usepackage{color}
\usepackage[colorlinks, linkcolor=blue, citecolor=blue]{hyperref}

\usepackage{sidecap}
\usepackage{float}
\usepackage{extarrows}
\usepackage{booktabs}
\usepackage{verbatim}
\usepackage[usenames,dvipsnames]{xcolor}

\newtheorem{theorem}{Theorem}[section]
\newtheorem{lemma}[theorem]{Lemma}

\theoremstyle{definition}
\newtheorem{remark}[theorem]{Remark}

\def\XXint#1#2#3{{\setbox0=\hbox{$#1{#2#3}{\int}$}
         \vcenter{\hbox{$#2#3$}}\kern-.5\wd0}}

\numberwithin{equation}{section}

\allowdisplaybreaks[4]

\begin{document}

\title{Sharp  blow-up result for the intercritical   Inhomogeneous  NLS equation }
\author{  Yuan  Li}
%\thanks{}
\date{}
\maketitle

\begin{abstract}
In this paper, we consider the intercritical inhomogeneous nonlinear Schr\"odinger equation.
For the radial symmetry initial data,   we construct the ring blow-up solutions and obtain blow-up speed. This result implies that the upper bound on the blow-up speed given by Cardoso and Farah [J. Funct. Anal.,281(8) No. 109134, (2021)] is sharp.
%In particular, our main result (and its proof) can be seen as a model scenario of ring blowup for mass super-critical nonlinear NLS with singular nonlinearity.

\noindent {\bf Keywords:} Inhomogeneous NLS equation; Intercritical regime;  Ring blow-up solution; Sharp upper bound
\end{abstract}

\section{Introduction and Main Result}
\noindent
We consider the  following  inhomogeneous nonlinear Schr\"odinger equation
\begin{equation}\label{equ:INLS-1}
\begin{cases}
     iu_t+\Delta u+|x|^{-\sigma}|u|^{p-1}u=0,~~(t,x)\in\mathbb{R}\times\mathbb{R}^N,\\
     u(t_0,x)=u_0(x),
     \tag{INLS}
    \end{cases}
\end{equation}
where $u(t,x)$ is a complex-valued function, $0<\sigma<\min\{\frac{N}{2},2\}$  and  $t_0\in\mathbb{R}$. This kind of problem arises naturally in nonlinear optics for the propagation of laser beams. We refer the reader to Gill \cite{G2000P}, Liu and Tripathi \cite{LT1994P} for more details.

Let us review some well-known results about \eqref{equ:INLS-1}. From Genoud and Stuart \cite{GS2008DCDS}, also see Guzman \cite{G2017NARWA}, given $u_0\in H^1(\mathbb{R}^N)$, there exists a unique solution $u\in C([t_0,T),H^1(\mathbb{R}^N))$ to \eqref{equ:INLS-1} and there holds the blowup alternative:
\begin{align*}
    T<+\infty,~~\text{implies}~~\lim_{t\to\infty}\|u(t)\|_{H^1}=+\infty.
\end{align*}
Moreover, the $H^1$ flow admits the conservation laws:

Mass:
\begin{align*}
    M(u(t))=\int|u(t,x)|^2dx=M(u_0);
\end{align*}

Energy:
\begin{align*}
    E_{\sigma}(u(t))=\frac{1}{2}\int|\nabla u(t,x)|^2dx-\frac{1}{p+1}\int|x|^{-\sigma}|u(t,x)|^{p+1}dx=E_{\sigma}(u_0).
\end{align*}
In addition, equation \eqref{equ:INLS-1} has the following symmetries
\[
u(t,x)\mapsto \lambda^{\frac{2-\sigma}{p-1}}u(\lambda^2t,\lambda x)e^{i\gamma},~~\lambda>0,~\gamma\in\mathbb{R}.
\]
It leaves invariant the norm in the homogeneous Sobolev space $\dot{H}^{s_c}$, where $s_c=\frac{N}{2}-\frac{2-\sigma}{p-1}$. If $s_c=0$, the problem is mass-critical, if $s_c=1$ it is energy-critical and  if $0<s_c<1$, it is mass super-critical and energy subcritical or just intercritical. The existence of solution with finite maximal time of existence is already known in $H^1$ for the equation \eqref{equ:INLS-1}. As for the classical NLS equation this is a consequence of the following Virial identity
\begin{align*}
    \frac{d^2}{dt^2}\int|x|^2|u(t,x)|^2=8((p-1)s_c+2)E(u_0)-4(p-1)s_c\|\nabla u(t)\|_{L^2}^2,
\end{align*}
satisfied by solutions to \eqref{equ:INLS-1} with initial data
\[u_0\in\Sigma=H^1\cap\{xu\in L^2\}.\]
This was obtained by Farah \cite{F2016JEE} and Dinh \cite{D2018NA}. From above identity, the solution blows
 up in finite time if $E_{\sigma}(u_0)<0$. In particular, Cardoso and Farah \cite{CF2021JFA,CF2023MZ} obtained that if the initial data $u_0\in H^1$ and the maximal time of existence $T^*>0$ for the corresponding solution $u\in ([0,T),H^1)$ of \eqref{equ:INLS-1} is finite, then the following space-time upper bound holds
 \begin{align}\label{upperbd-1}
     \int_t^{T^*}(T^*-\tau)\|\nabla u(\tau)\|_{L^2}^2d\tau\leq C(u_0)(T^*-t)^{\frac{2(5-p)}{(p-1)(N-2)+2\sigma+(5-p)}},
 \end{align}
for $t$ close enough to $T^*$ ($N\geq3$ and $0<\sigma<1$). For more details about the existence of  blow-up solution and blow-up rate, one can see \cite{CF2023MZ,Z2014JMAA,PZ2018MMAS,ADF2021CPDE} and the references therein.

If $\sigma=0$, then \eqref{equ:INLS-1} is the classical  NLS equation
\begin{equation}\label{equ:NLS-1}
\begin{cases}
     iu_t+\Delta u+|u|^{p-1}u=0,~~(t,x)\in\mathbb{R}\times\mathbb{R}^N,\\
     u(t_0,x)=u_0(x).
    \end{cases}
\end{equation}
 From Ginbre and Velo \cite{GV1979JFA}, given $u_0\in H^1$, there exists a unique solution $u\in C([t_0,T),H^1)$ and there holds the blow-up alternative. Moreover, the   $H^1$ flow satisfies the conservation laws of mass, energy, momentum. Further more, a group of $H^1$ symmetries leaves the flow invariant: if $u(t,x)$ solves \eqref{equ:NLS-1} then for any $(\lambda_0,\tau_0,x_0,\gamma_0)\in\mathbb{R}^+\times\mathbb{R}\times\mathbb{R}^N\times\mathbb{R}$, so does
\begin{align*}
    v(t,x)\mapsto\lambda_0^{\frac{2}{p-1}}u(\lambda_0^2t+\tau_0,\lambda_0x+x_0)e^{i\gamma_0},
\end{align*}
and the Galilean drift
\begin{align*}
    v(t,x)\mapsto u(t,x-\beta_0t)e^{i\frac{\beta_0}{2}\cdot\left(x-\frac{\beta_0}{2}t\right)},~~\beta_0\in\mathbb{R}^N.
\end{align*}
Notice that by the scaling invariant homogeneous Sobolev space $\dot{H}^{s_c^0}$, we see that $s_c^0=\frac{N}{2}-\frac{2}{p-1}>0$. If $0<s_c^0<1$, this means that the problem \eqref{equ:NLS-1} is mass super-critical and energy-subcritical. If the initial data $u_0\in\sum=H^1\cap\{xu\in L^2\}$, then from \cite{Book2003NLS}, the following Virial identity holds
\begin{align}\label{virial:intro-1}
    \frac{d^2}{dt^2}\int|x|^2|u(t,x)|^2dx=4N(p-1)E_0(u_0)-\frac{16s_c^0}{N-2s_c^0}\int|\nabla u|^2\leq 4N(p-1)E_{0}(u_0).
\end{align}
This implies that if the initial energy $ E_0(u_0)<0$, then the solution of problem \eqref{equ:NLS-1}  blows up in finite time. Also, \eqref{virial:intro-1} implies that there exists a universal upper bound on blow-up rate (see \cite{Book2003NLS}):
\begin{align*}
    \int_0^T(T-t)\|\nabla u(t)\|_{L^2}^2dt<+\infty.
\end{align*}
Furthermore, Merle, Rapha\"el and Szeftel \cite{MRS2014Duke} obtained the  upper bound on the blow-up speed for the radial data ($N\geq2,~p\leq5$),
\begin{align}\label{blowup:upper-1}
    \int_t^T(T-\tau)\|\nabla u(\tau)\|_{L^2}^2d\tau\leq C(T-t)^{\frac{2(5-p)}{(p-1)(N-1)+5-p}}.
\end{align}
They also constructed  the collapsing ring blow-up solution ($N\geq2$, $1+\frac{2}{N}<p<\min\left\{\frac{N+2}{N-2},5\right\}$),
\[u(t,x)\sim\frac{1}{\lambda(t)^{\frac{2}{p-1}}}[Qe^{-i\beta_0y}]\left(\frac{r-\alpha(t)}{\lambda(t)}\right)e^{i\gamma},\]
where
\begin{align*}
    \lambda(t)\sim(T-t)^{1+\frac{5-p}{(p-1)(N-1)}},~~~\alpha(t)\sim(T-t)^{\frac{(5-p)}{(p-1)(N-1)+5-p}},
\end{align*}
and $Q$ (see \cite{K1989ARMA,GNN1979CMP}) is the one-dimensional unique positive ground state of \begin{align}\label{equ:elliptic-1}
    -\Delta Q+Q+Q^p=0,
\end{align}
and proved that the upper bound \eqref{blowup:upper-1} is sharp.

Inspired by \cite{MRS2014Duke}, in this paper, we aim to construct the ring blowup solution to the inhomogeneous NLS equation
\begin{align}\label{equ:3-1}
    iu_t+\Delta u+|x|^{-\sigma}|u|^{2}u=0,~~~(t,x)\in\mathbb{R}\times\mathbb{R}^3,
\end{align}
and prove the upper bound \eqref{upperbd-1} is sharp.

 Throughout this paper, we define
\begin{align}\label{def:alpha0-1}
    a_{\infty}=\frac{1}{2(\sigma+1)},~~
    \beta_\infty=\sqrt{\frac{1}{3(1+\sigma)}}.
\end{align}
Now we state our main result.

\begin{theorem}\label{thm-1}

Let $0<\sigma<1$ and  $Q$ be the one-dimensional mass subcritical ground state solution to \eqref{equ:elliptic-1}($p=3$). Then there exist a time $t_0<0$ and a solution $u\in C([t_0,0),H^1)$ of \eqref{equ:3-1} with radial symmetry which blows up at time $T=0$. More precisely, it holds that
\begin{align*}
    u(t,x)-(\alpha(t))^{\frac{\sigma}{2}}\frac{1}{\lambda(t)}\left[Qe^{-i\beta_{\infty}y}\right]\left(\frac{r-\alpha(t)}{\lambda(t)}\right)e^{i\gamma(t)}\to0~~\text{in}~~L^2(\mathbb{R}^3)~~\text{as}~~t\to0^-,
\end{align*}
where $y=\frac{r-\alpha(t)}{\lambda(t)}$ and
\[\lambda(t)\sim|t|^{\frac{1}{1+\frac{ a_{\infty}}{1-\sigma a_{\infty}}}},~~\alpha(t)\sim|t|^{\frac{\frac{ a_{\infty}}{1-\sigma a_{\infty}}}{1+\frac{ a_{\infty}}{1-\sigma a_{\infty}}}},~~\gamma(t)\sim|t|^{-\frac{1-(\sigma+1) a_{\infty}}{1-(\sigma-1) a_{\infty}}},~~\text{as}~~t\to0^-.\]
Moreover, the blowup speed is given by
\begin{align}\label{blow:rate-1}
    \|\nabla u(t)\|_{L^2}\sim\frac{1}{|T^*-t|^{\frac{2+\sigma}{3+\sigma}}}~~\text{as}~~t\to T^*.
\end{align}

\end{theorem}

{\bf Comments:}

1. Extension. Similar result can be addressed for the problem \eqref{equ:INLS-1} with $N\geq2$, $\sigma\in(0,1)$ and $1+\frac{4-2\sigma}{N}<p<\min\left\{\frac{N+2-2\sigma}{N-2},5\right\}$.  We claim that blowup solution have the following form
\begin{align*}
    u(t,x)-(\alpha(t))^{\frac{\sigma}{2}}\frac{1}{\lambda^{\frac{2}{p-1}}(t)}\left[Qe^{-i\beta^1_{\infty}y}\right]\left(\frac{r-\alpha(t)}{\lambda(t)}\right)e^{i\gamma(t)}\to0~~\text{in}~~L^2(\mathbb{R}^N)~~\text{as}~~t\to0^-,
\end{align*}
where $y=\frac{r-\alpha(t)}{\lambda(t)}$ and
\[\lambda(t)\sim|t|^{\frac{1}{1+\frac{ a^1_{\infty}}{1-\sigma a^1_{\infty}}}},~~\alpha(t)\sim|t|^{\frac{\frac{ a^1_{\infty}}{1-\sigma a^1_{\infty}}}{1+\frac{ a^1_{\infty}}{1-\sigma a^1_{\infty}}}},~~\gamma(t)\sim|t|^{-\frac{1-(\sigma+1) a^1_{\infty}}{1-(\sigma-1) a^1_{\infty}}},~~\text{as}~~t\to0^-.\]
Here
\begin{align*}
    a^1_{\infty}=\frac{5-p}{(p-1)(N-1)+4\sigma},~~\beta^1_\infty=\sqrt{\frac{2(p-1)(5-p)}{(p+3)((p-1)(N-1)+4\sigma)}}.
\end{align*}
The corresponding blowup speed speed is
\begin{align*}
    \|\nabla u(t)\|_{L^2}\sim\frac{1}{|T^*-t|^{\frac{1}{1+\frac{5-p}{(p-1)(N-1+\sigma)}}}}~~\text{as}~~t\to T^*.
\end{align*}
The proof carries over verbatim except for some technicalities when the nonlinearity $|x|^{-\sigma}|u|^{p-1}u$ fails to be smooth.

2. Sharp upper bound on the blowup speed. The blowup rate \eqref{blow:rate-1} of ring solutions saturates the upper bound \eqref{upperbd-1}, which is therefore optimal in the radial setting.

3. Although in \cite{RS2011JAMS,BCD2011CPDE}, they considered the minimal mass blowup solution of  inhomogeneous  mass critical NLS equation, where the inhomogeneous factor is smooth, bounded and integrable. However, in our case, the inhomogeneous term $\frac{1}{|x|^{\sigma}}$ is singular, so we need to use  some new techniques to handle it. Fortunately, our problem \eqref{equ:INLS-1}  has the scaling symmetry compared to the equations in \cite{RS2011JAMS,BCD2011CPDE}. This is very useful for us to construct the approximate profile.

4. Due to the inhomogeneous factor, the parameters are different from the classical NLS equation, see Lemma \ref{lemma:app:pro-1} and \ref{lemma:initial-1}, so we need to be careful to obtain the sharp result.

\begin{remark}
    Notice that in Theorem \ref{thm-1},
    $Q$  is the ground state of \eqref{equ:elliptic-1} instead of the ground state of
    \begin{align}\label{equ:e2-1}
        -\Delta u+u+\frac{1}{|x|^\sigma}u^p=0.
    \end{align}
One reason is that the inhomogeneous factor $\frac{1}{|x|^\sigma}$ is a singular term, the ground state of equation \eqref{equ:e2-1} is  undefined at origin. Another important reason is that it is very difficult to construct approximate solutions by using the ground state of \eqref{equ:e2-1}. In our argument, we treat $\frac{1}{|x|^{\sigma}}$ as a weight function, so we can use the ground state  $Q$ to construct the approximate solution of the equation \eqref{equ:3-1}. Actually, it is a naturally choice to construct the approximate solution  by using  the ground state solution $Q$ of equation \eqref{equ:elliptic-1}(see Section 2).

\end{remark}

{\bf Notations}

For the positive $a$ and $b$, the notation $a\lesssim b$ means that $a\leq cb$ holds for a universal constant $c>0$, $a\gtrsim b$ means that $a\geq cb$ holds for a universal constant $c>0$, $a\sim b$ means that $a\lesssim b$ and $a\gtrsim b$. Let $(\cdot,\cdot)$ denote the scalar product on $L^2$,
\[(f,g)=\int f(x)\bar{g}(x)dx,\]
for $f,g$ two complex valued functions in $L^2$.

Denote the operator
\begin{align*}
    \Lambda f=f+x\cdot\nabla f.
\end{align*}
Let $L=(L_+,L_-)$ be the linearized operator around $Q$,
\begin{equation*}
\begin{aligned}
L_+=-\partial_y^2+1-3Q^2,~~
    L_-=-\partial_y^2+1-Q^2,
\end{aligned}
\end{equation*}
where $Q$ be the one-dimensional mass sub-critical ground state solution to \eqref{equ:elliptic-1}.

The rest of the paper is as follows. In section 2, we construct an approximate solution $Q_{\mathcal{P}}$ of the renormalized solution. In section 3, we decompose the solution and estimate the modulation parameters. In section 4, we establish a refined energy/virial type estimate. In section 5, we apply the energy estimate to establish a bootstrap argument that will be needed in the construction of solutions that ring blow up. In section 6, we give the proof of the Theorem \ref{thm-1}. The finally section is the appendix.

\section{The approximate solution}
In this section, we aim to construct the approximate solution at any order.
Let us consider the general modulated
\[u(t,x)=\alpha^{\frac{\sigma}{2}}(t)\frac{1}{\lambda(t)}v\left(s,\frac{r-\alpha(t)}{\lambda(t)}\right)e^{i\gamma(t)},~~~\frac{ds}{dt}=\frac{1}{\lambda^2}.\]
Injecting this into the problem \eqref{equ:3-1} yields:
\begin{align}\label{equ:scaling-1}
    i\partial_sv+v_{yy}+\frac{\lambda}{\alpha}\frac{2}{1+\frac{\lambda}{\alpha}y}v_y-v+\frac{1}{\left|\frac{\lambda}{\alpha} y+1\right|^\sigma}|v|^2v=i\frac{\lambda_s}{\lambda}\Lambda v+i\frac{\alpha_s}{\lambda}\left(v_y-\frac{\sigma}{2}\frac{\lambda}{\alpha}v\right)+\Tilde{\gamma}_sv,
\end{align}
where $\Tilde{\gamma}_s=\gamma_s-1$.

Now we define
\begin{align}\label{def:a-1}
    a=\frac{2\beta}{a_{\infty}}\frac{\lambda}{\alpha}.
\end{align}
Then, \eqref{equ:scaling-1} is equivalent to
\begin{align*}
     i\partial_sv+v_{yy}+\frac{2}{1+\frac{a_{\infty}a}{2\beta}y}\frac{a_{\infty}a}{2\beta}v_y-v+\frac{1}{\left|1+\frac{a_{\infty}a}{2\beta}y\right|^\sigma}|v|^2v=i\frac{\lambda_s}{\lambda}\Lambda v+i\frac{\alpha_s}{\lambda}\left(v_y-\frac{\sigma a_{\infty}a}{4\beta}v\right)+\Tilde{\gamma}_sv,
\end{align*}
Let
\begin{align*}
    w(s,y)=v(s,y)e^{i\beta y},
\end{align*}
which satisfies
\begin{align}\label{equ:w-1}
    &i\partial_sw+w_{yy}-w+\frac{1}{\left|1+\frac{ a_{\infty}a}{2\beta}y\right|^\sigma}|w|^2w+\frac{a_{\infty}a}{\beta}\frac{1}{1+\frac{a_{\infty}a}{2\beta}y}(w_y-i\beta w)+a(i\Lambda w+\beta yw)-\frac{i\sigma a_{\infty}}{2}aw\notag\\
    =&-\beta_syw+\left(\frac{\lambda_s}{\lambda}+a\right)(i\Lambda w+\beta yw)+\left(\frac{\alpha_s}{\lambda}+2\beta\right)\left(iw_y+\beta w-\frac{i\sigma a_{\infty}a}{4\beta}w\right)+(\Tilde{\gamma}_s-\beta^2)w.
\end{align}
Let
\[\beta=\beta_{\infty}+\Tilde{\beta}.\]
We look for an approximate solution to \eqref{equ:scaling-1} of the form
\[v(s,y)=Q_{\mathcal{P}}(y),~~~\text{where}~~~\mathcal{P}=(a,\Tilde{\beta}).\]
Following the
slow modulated ansatz strategy developed in \cite{RS2011JAMS,MRS2014Duke,L2022}, we freeze the modulation
\[
\frac{\lambda_s}{\lambda}=-a+\mathbf{A}_1(\mathcal{P}),~~\frac{\alpha_s}{\lambda}=-2\beta,~~\Tilde{\gamma}_s=\beta^2,~~~\beta_s=\mathbf{A}_2(\mathcal{P}),\]
where $\mathbf{A}_1$ and $\mathbf{A}_2$ are polynomials in $\mathcal{P}$, which will be chosen later to ensure suitable solvability conditions. From the definition of $a$ (see \eqref{def:a-1}), we have the following relation
\begin{align*}
    &a_s+(1-a_{\infty})a^2-\frac{a}{\beta}\mathbf{A}_2-a\mathbf{A}_1\notag\\
    =&\frac{a}{\beta}(\Tilde{\beta}_s-\mathbf{A}_2)+a\left(\frac{\lambda_s}{\lambda}+a-\mathbf{A}_1\right)-\frac{a_{\infty}}{2\beta}a^2\left(\frac{\alpha_s}{\lambda}+2\beta\right).
\end{align*}
From the above, we have
\begin{align}\label{def:error-1}
    &i\left(-(1-a_{\infty})a^2+\frac{a}{\beta}\mathbf{A}_2+a\mathbf{A}_1\right)\partial_aQ_{\mathcal{P}}+i\mathbf{A}_2\partial_{\Tilde{\beta}} Q_{\mathcal{P}}\notag\\
    &-(1+\beta^2)Q_{\mathcal{P}}+i(a-\mathbf{A}_1)(1+y\partial_y)Q_{\mathcal{P}}+2i\beta\partial_yQ_{\mathcal{P}}-i\frac{\sigma a_{\infty}}{2}a Q_{\mathcal{P}}\notag\\
    &+\partial_y^2Q_{\mathcal{P}}+\frac{a_{\infty}a}{\beta}\frac{1}{1+\frac{a_{\infty}a}{2\beta}y}\partial_yQ_{\mathcal{P}}+\frac{1}{\left|1+\frac{a_{\infty}a}{2\beta}y\right|^\sigma}|Q_{\mathcal{P}}|^2Q_{\mathcal{P}}
    =-\Psi_{\mathcal{P}}.
\end{align}
Let
\begin{align*}
    Q_{\mathcal{P}}(y)=P_{\mathcal{P}}e^{-i\beta y-ia\frac{y^2}{4}},
\end{align*}
so that \eqref{def:error-1}
\begin{align}\label{equ:app-1}
&i\left(-(1-a_{\infty})a^2+\frac{a}{\beta}\mathbf{A}_2+a\mathbf{A}_1\right)\partial_aP_{\mathcal{P}}+i\mathbf{A}_2\partial_{\Tilde{\beta}} P_{\mathcal{P}}\notag\\
    &+\partial_y^2P_{\mathcal{P}}-P_{\mathcal{P}}+\frac{a_{\infty}a}{\beta}\frac{1}{1+\frac{a_{\infty}a}{2\beta}y}\partial_yP_{\mathcal{P}}+\frac{1}{\left|1+\frac{a_{\infty}a}{2\beta}\right|^\sigma}|P_{\mathcal{P}}|^2P_{\mathcal{P}}\notag\\
    &-i\mathbf{A}_1(1+y\partial_y) P_{\mathcal{P}}-\mathbf{A}_1\left(\beta y+\frac{ay^2}{2}\right)P_{\mathcal{P}}+\mathbf{A}_2yP_{\mathcal{P}}\notag\\
    &+\left[a\beta y+\left(a_{\infty}a^2+\frac{a}{\beta}\mathbf{A}_2+a\mathbf{A}_1\right)\frac{y^2}{4}-i\left(\frac{1}{1+\frac{a_{\infty}a}{2\beta}y}\frac{a_{\infty}a^2}{2\beta}(1-a_{\infty})y\right)-\frac{i\sigma a_{\infty}}{2}a\right]P_{\mathcal{P}}\notag\\
    =&-\Psi_{\mathcal{P}}e^{i\beta y+ia\frac{y^2}{4}}.
\end{align}
Next Lemma we will construct the  approximate solution to \eqref{equ:app-1}.
\begin{lemma}\label{lemma:app:pro-1}
(Approximate Profile). Let an integer $l>5$, then there exist polynomials $\mathbf{A}_1$ and $\mathbf{A}_2$ of the following form
\begin{equation}\label{def:a1a2-1}
    \begin{aligned}
    \mathbf{A}_1(\mathcal{P})&=\sigma a_{\infty} a+\sum_{2\leq j+k\leq l-1}c_{1,j,k}a^j\Tilde{\beta}^k,\\
    \mathbf{A}_2(\mathcal{P})&
    =\sigma a_{\infty}\beta_{\infty}a+\sum_{2\leq j+k\leq l-1}c_{2,j,k}a^j\Tilde{\beta}^k,
    \end{aligned}
\end{equation}
where $c_{1,j,k}$ and $c_{2,j,k}$ are the constants and smooth well-localized solutions $(T_{j,k},S_{j,k})$, such that
\begin{align}\label{app:sol:form-1}
    P_{\mathcal{P}}=Q+\sum_{1\leq j+k\leq l-1}a^j\Tilde{\beta}^k(T_{j,k}+iS_{j,k}),
\end{align}
where $Q$ is the one dimensional ground state solution to \eqref{equ:elliptic-1}$(p=3)$, is a approximate solution to \eqref{equ:app-1} with $\Psi_{\mathcal{P}}$ smooth and well-localized in $y$ satisfying the following decay estimate
\begin{align*}
    \Psi_{\mathcal{P}}=\mathcal{O}(a^l|y|^{c_l}e^{-|y|}).
\end{align*}
Moreover, the approximate profile holds the decay estimate
\begin{align}\label{decay:P-1}
    |P_{\mathcal{P}}|\lesssim(1+|y|^{2l})e^{-|y|}.
\end{align}
\end{lemma}
\begin{proof}
From \cite{W1985SIAM,CGNT2007SIAM}, the kernel of $L=(L_+,L_-)$
is explicit
\begin{align}\label{kernel-1}
    \ker L_+=span\{Q^\prime\}.~~\ker L_-=span\{Q\}.
\end{align}
It follows from the kernel properties \eqref{kernel-1} of $L_+$ and $L_-$ and well-known properties of the Helmholtz kernel that
\begin{equation*}
    \begin{aligned}
    \forall~g\in H^1(\mathbb{R}),~~(g,Q^\prime)=0,~~\exists~f_+\in H^1(\mathbb{R}),~s.~t.~L_+f_+=g,\\
    \forall~g\in H^1(\mathbb{R}),~~(g,Q)=0,~~\exists~f_-\in H^1(\mathbb{R}),~s.~t.~L_+f_-=g.
    \end{aligned}
\end{equation*}
We also give the following Pohozaev indentities
\begin{equation}\label{Pohozaev:iden-1}
    \begin{aligned}
     &-2\int (Q^\prime)^2+2\int Q^2=\int Q^4,\\
     &\int (Q^\prime)^2+\int Q^2=\int Q^4.
    \end{aligned}
\end{equation}
From the above identities \eqref{Pohozaev:iden-1}, we have
\begin{align}\label{identity-1}
    3\int (Q^\prime)^2=\int Q^2=\frac{3}{4}\int Q^4.
\end{align}
We recall that $L=(L_+,L_-)$ has a generalized nullspace characterized by the following algebraic identities generated by the symmetry group:
\begin{align*}
&L_-Q=0,~~L_-{yQ}=-2\nabla Q,\\
&L_+Q^\prime=0,~~L_+\Lambda Q=-2Q.
\end{align*}

The proof proceeds by injecting the expansion \eqref{app:sol:form-1} into \eqref{equ:app-1}, identifying the terms with the same homogeneity, and inverting the corresponding operator. Notice that we have the following Taylor expansion
\begin{align*}
    (1+x)^b=1+bx+\frac{b(b-1)}{2!}x^2+\frac{b(b-1)(b-2)}{3!}x^3+\ldots.
\end{align*}

Now we divide the proof into the following steps.

{\bf Step 1:} General case. Let $j+k\geq1$. Assume that $T_{p,q}$, $S_{p,q}$, $c_{1,p,q}$ and $c_{2,p,q}$ for $q+p\leq j+k-1$ have been constructed. Then, identifying the terms homogeneous of order $(j,k)$ in \eqref{equ:app-1} yields a linear system of the following type
\begin{align}\label{system:linear-1}
    \begin{cases}
    L_+T_{j,k}=h_{1,j,k}-c_{1,j,k}\beta_{\infty}yQ+c_{2,j,k}yQ,\\
    L_-S_{j,k}=h_{2,j,k}-c_{1,j,k}\Lambda Q,
    \end{cases}
\end{align}
where $h_{1,j,k}$ and $h_{2,j,k}$ can be computed explicitly and only depend on $T_{p,q}$, $S_{p,q}$, $c_{1,p,q}$ and $c_{2,p,q}$ for $q+p\leq j+k-1$. According to \eqref{kernel-1}, the invertibility of \eqref{system:linear-1} must be satisfy the orthogonality condition
\begin{align}\label{condition:solve-1}
\begin{cases}
\left(h_{1,j,k}-c_{1,j,k}\beta_{\infty}yQ+c_{2,j,k}yQ,Q^{\prime}\right)=0,\\
\left(h_{2,j,k}-c_{1,j,k}\Lambda Q,Q\right)=0.
\end{cases}
\end{align}
Now we claim: For all $1\leq j+k\leq l-1$, let
\begin{align}\label{def:c1-1}
    &c_{1,j,k}=\frac{1}{(Q,\Lambda Q)}(h_{2,j,k},Q),\\\label{def:c2-1}
    &c_{2,j,k}=\frac{2}{\|Q\|_{L^2}^2}(h_{1,j,k},Q^{\prime})+\frac{\beta_{\infty}}{(Q,\Lambda Q)}(h_{2,j,k},Q).
\end{align}
Then, there exist $(T_{j,k},S_{j,k})$ solution of \eqref{equ:app-1} for all $1\leq j+k\leq l-1$. Furthermore, $T_{j,k}$ and $S_{j,k}$ are smooth and decay as
\begin{align*}
    T_{j,k}, S_{j,k}=\mathcal{O}(|y|^{2(j+k)}e^{-|y|}).
\end{align*}
In fact, to be able to solve for $(T_{j,k},S_{j,k})$. If we choosing $c_{1,j,k}$ and $c_{2,j,k}$ as in \eqref{def:c1-1} and \eqref{def:c2-1}, respectively, we may solve for $(T_{j,k},S_{j,k})$ solution of \eqref{system:linear-1}.

Next, we investigate the smoothness and decay properties of $(T_{j,k},S_{j,k})$. Identifying the terms homogeneous of order $j+k$ in \eqref{equ:app-1}, we have for $h_{1,j,k}$ and $h_{2,j,k}$ defined in \eqref{system:linear-1}
\begin{align*}
    \begin{cases}
    h_{1,j,k}=\sum\limits_{p+q\leq j+k-1}\left(a_{1,p,q}y^{j+k-p-q}T_{p,q}+a_{2,p,q}y^{j+k-p-q}S_{p,q}+a_{3,p,q}y^{j+k-p-q}T^{\prime}_{p,q}\right)+NL_j^{(1)},\\
    h_{2,j,k}=\sum\limits_{p+q\leq j+k-1}\left(a_{4,p,q}y^{j+k-p-q}T_{p,q}+a_{5,p,q}y^{j+k-p-q}S_{p,q}+a_{6,p,q}y^{j+k-p-q}T^{\prime}_{p,q}\right)+NL_j^{(2)},
    \end{cases}
\end{align*}
where $T_{0,0}=Q$, $a_{n,p,q}$ are real numbers which may be explicitly computed, and where $NL_j^{(1)}$ and $NL_j^{(2)}$ are the contributions coming from the Taylor expansion of the  term $\frac{1}{|1+\frac{a_{\infty}a}{2\beta}y|^{\sigma}}|P_{\mathcal{P}}|^2P_{\mathcal{P}}$.
Then ,by the similar argument as \cite{MRS2014Duke}, we can obtain the smoothness and decay properties of $T_{j,k}$ and $S_{j,k}$. Here we omit it.

{\bf Step 2:} Computation of $c_{1,1,0}$ and $c_{2,1,0}$. From the terms homogeneous of order $(1,0)$ in \eqref{equ:app-1}, we get
\begin{align*}
    \begin{cases}
    L_+T_{1,0}=\frac{a_{\infty}}{\beta_{\infty}}Q^\prime+\beta_{\infty} yQ-c_{1,1,0}\beta_{\infty}yQ+c_{2,1,0}yQ+\frac{\sigma a_{\infty}}{2\beta_{\infty}}Q^3,\\
    L_-S_{1,0}=-c_{1,1,0}\Lambda Q-\frac{\sigma a_{\infty}}{2}Q.
    \end{cases}
\end{align*}
From \eqref{kernel-1} and \eqref{condition:solve-1},  the solvability
conditions for $T_{1,0}$ and $S_{1,0}$ are equivalent to
\begin{align*}
    \begin{cases}
    \left(\frac{a_{\infty}}{\beta_{\infty}}Q^\prime+\beta_{\infty} yQ-c_{1,1,0}\beta_{\infty}yQ+c_{2,1,0}yQ+\frac{\sigma a_{\infty}}{2\beta_{\infty}}Q^3,Q^\prime\right)=0,\\
    \left(-c_{1,1,0}\Lambda Q-\frac{\sigma a_{\infty}}{2}Q,Q\right)=0.
    \end{cases}
\end{align*}
By the second equation, we can obtain \begin{align*}
    c_{1,1,0}=\frac{\sigma a_{\infty}}{2(\Lambda Q,Q)}\|Q\|_{L^2}^2=\sigma a_{\infty}.
\end{align*}
Notice that
\begin{align*}
    \left(\frac{a_{\infty}}{\beta_{\infty}}Q^\prime+\beta_{\infty} yQ,Q^\prime\right)=\frac{a_{\infty}}{\beta_{\infty}}\int (Q^\prime)^2-\frac{\beta_{\infty}}{2}\int Q^2=\left(\frac{a_{\infty}}{\beta_{\infty}}-\frac{3\beta_{\infty}}{2}\right)\int (Q^\prime)^2=0,
\end{align*}
where we used the equality \eqref{identity-1} and the definition of $a_{\infty}$ and $\beta_{\infty}$ (see \eqref{def:alpha0-1}). Hence, we can obtain that
\begin{align*}
    c_{2,1,0}=\beta_{\infty}c_{1,1,0},
\end{align*}
since $Q$ is even function and $Q^\prime$ is odd.

{\bf Step 3:} Computation of $c_{1,0,1}$ and $c_{2,0,1}$. From the terms homogeneous of order $(0,1)$ in \eqref{equ:app-1} and get
\begin{align*}
    \begin{cases}
    L_+T_{0,1}=-c_{1,0,1}\beta_{\infty}yQ+c_{2,0,1}yQ,\\
    L_-S_{0,1}=-c_{1,0,1}\Lambda Q.
    \end{cases}
\end{align*}
By \eqref{condition:solve-1}, we can obtain
\begin{align*}
    c_{1,0,1}=0=c_{2,0,1}.
\end{align*}

{\bf Step 4:} Computation of $c_{1,2,0}$ and $c_{2,2,0}$. From the terms homogeneous of order $(2,0)$ in \eqref{equ:app-1} and get
\begin{align*}
    \begin{cases}
    L_+T_{2,0}=&\left(-(1-a_{\infty})+c_{1,1,0}+\frac{c_{2,1,0}}{\beta_{\infty}}\right)S_{1,0}+\frac{a_{\infty}}{\beta_{\infty}}T^{\prime}_{1,0}-\frac{a_{\infty}^2}{2\beta_{\infty}^2}yQ^\prime+3QT_{1,0}^2+QS_{1,0}^2+\frac{\sigma a_{\infty}}{2\beta_{\infty}}3Q^2T_{1,0}\\
    &+\frac{\sigma(\sigma-1)a_{\infty}^2}{8\beta_{\infty}^2}Q^3+\beta_{\infty} yT_{1,0}+\left(a_{\infty}+\frac{c_{2,1,0}}{\beta_{\infty}}+c_{1,1,0}\right)\frac{1}{4}y^2Q+\frac{\sigma a_{\infty}}{2}S_{1,0}\\
    &-c_{1,2,0}\beta_{\infty}yQ+c_{2,2,0}yQ+\frac{\sigma a_{\infty}}{2}S_{1,0},\\
    L_-S_{1,0}=&\left(-(1-a_{\infty})+c_{1,1,0}+\frac{c_{2,1,0}}{\beta_{\infty}}\right)T_{1,0}+\frac{a_{\infty}}{\beta_{\infty}}S^{\prime}_{1,0}+2QT_{1,0}S_{1,0}+\frac{\sigma a_{\infty}}{2\beta_{\infty}}Q^2S_{1,0}+\beta_{\infty} yS_{1,0}\\
    &-\frac{a_{\infty}}{2\beta_\infty}(1-\frac{2\sigma a_{\infty}}{2}T_{1,0}-a_{\infty})yQ-c_{1,2,0}\Lambda Q.
    \end{cases}
\end{align*}
Notice that $T_{1,0}$ is an odd function, while $S_{1,0}$ is an even function. From the solvability conditions for $T_{2,0}$ and $S_{2,0}$, we can obtain
\begin{align*}
    c_{1,2,0}\neq0,~~c_{2,2,0}\neq0.
\end{align*}

{\bf Step 5:} Conclusion. The error term $\Psi_{\mathcal{P}}$ consists of polynomial in $(T_{j,k},S_{j,k})_{j+k\leq l-1}$ with lower-order $l$, the error between the Taylor expansion of the potential terms $\frac{a_{\infty}a}{2\beta}\frac{2}{1+\frac{a_{\infty}a}{2\beta}y}$  and $\frac{2}{1+\frac{a_{\infty}a}{2\beta}y}$ in \eqref{equ:app-1}, and between the  Taylor expansion of the nonlinear term $\frac{1}{\left|1+\frac{a_{\infty}a}{2\beta}\right|^{\sigma}}|P_{\mathcal{P}}|^2P_{\mathcal{P}}$ and $\frac{1}{|y|^b}$. Using the exponential decay of $P_{\mathcal{P}}$, we can easily estimate these terms.
\end{proof}

\section{Modulation Estimate}
In this section, we aim to give the parameters estimate that will be useful to prove the existence of the blowup solution.
Now we introduce a smooth cut-off function
\begin{align*}
    \xi(y)=\begin{cases}
         0~~&\text{for}~~y\leq-2,\\
         1~~&\text{for}~~y\geq-1.
    \end{cases}
\end{align*}
Let $\xi_a(y)=\xi(\sqrt{a}y)$ and define
\begin{align}\label{def:Qcut-1}
    \mathbf{Q}_{\mathcal{P}}(y)=\xi_a(y)P_{\mathcal{P}}(y)e^{-i\beta y-ia\frac{y^2}{4}}.
\end{align}

Given $C^1$ modulation parameters $(\lambda(t),\alpha(t),\gamma(t),\Tilde{\beta}(t))$ such that $0<a=\frac{2\beta}{a_{\infty}}\frac{\lambda(t)}{\alpha(t)}\ll1$, let
\begin{align}\label{def:Qtilde-1}
    \Tilde{Q}(t,x)=\alpha^{\frac{\sigma}{2}}(t)\frac{1}{\lambda(t)}\mathbf{Q}_{\mathcal{P}}\left(s,\frac{r-\alpha(t)}{\lambda(t)}\right)e^{i\gamma(t)}.
\end{align}
Then, by Lemma \ref{lemma:app:pro-1}, we can obtain that  $\Tilde{Q}$ is a smooth radially symmetric function which satisfies
\begin{align}\label{equ:Qtilde-1}
    i\partial_t\Tilde{Q}+\Delta\Tilde{Q}+|x|^{-\sigma}|\Tilde{Q}|^2\Tilde{Q}=\psi=\alpha^{\frac{\sigma}{2}}(t)\frac{1}{\lambda^{3}(t)}\Psi\left(t,\frac{r-\alpha(t)}{\lambda(t)}\right)e^{i\gamma(t)}
\end{align}
with
\begin{align}\label{def:Psi-1}
    \Psi=&-(\gamma_s-1-\beta^2)\mathbf{Q}_{\mathcal{P}}-i\left(\frac{\lambda_s}{\lambda}-a-\mathbf{A}_1\right)(\Lambda \mathbf{Q}_{\mathcal{P}}-a\partial_a\mathbf{Q}_{\mathcal{P}})\notag\\
    &-i\left(\frac{\alpha_s}{\lambda}+2\beta\right)\left(\partial_y\mathbf{Q}_{\mathcal{P}}+\frac{a_{\infty}}{2\beta}a^2\partial_a\mathbf{Q}_{\mathcal{P}}-\frac{\sigma a_{\infty}a}{4\beta}\mathbf{Q}_{\mathcal{P}}\right)\notag\\
    &+\left(a_s+(1-a_{\infty})a^2-\frac{a}{\beta}\mathbf{A}_2-a\mathbf{A}_1\right)\partial_a\mathbf{Q}_{\mathcal{P}}\notag\\
    &+i(\beta_s-\mathbf{A}_2)\left(\partial_{\Tilde{\beta}}\mathbf{Q}_{\mathcal{P}}+\frac{a}{\beta}\partial_a\mathbf{Q}_{\mathcal{P}}\right)\notag\\
    &+\mathcal{O}\left(\frac{e^{-|y|}}{a^{c_l}}\mathbf{1}_{y\sim\frac{1}{\sqrt{a}}}+|\mathcal{P}|^l\xi_a|y|^{c_l}e^{-|y|}\right).
\end{align}

\subsection{The approximation of the parameters}
In this subsection, we give the exact modulation equations formally predicted by the $Q_{\mathcal{P}}$ construction. Now we have the following lemma.
\begin{lemma}\label{lemma:initial-1}
There exist $t_0<0$ small enough and a solution $(\lambda,a,\Tilde{\beta}, \alpha,\gamma)$ to the dynamical system
\begin{align}\label{dynamical:sys-1}
    \begin{cases}
    \frac{\lambda_s}{\lambda}=-a+\mathbf{A}_1(\mathcal{P}),~~\mathcal{P}=(a,\Tilde{\beta}),\\
    \frac{\alpha_s}{\lambda}=-2\beta,~~\beta=\beta_{\infty}+\Tilde{\beta},\\
    \beta_s=\mathbf{A}_2,\\
    \gamma_s=1+\beta^2,\\
    \frac{ds}{dt}=\frac{1}{\lambda^2},\\
    a=\frac{2\beta}{a_{\infty}}\frac{\lambda}{\alpha},
    \end{cases}
\end{align}
which is defined on $[t_0,0)$. Moreover, this solution satisfies the following bounds,
\begin{align}\label{bound:a1-1}
&a(t)=\frac{1}{(1-(2\sigma+1) a_{\infty})}B_1^{\frac{1-(\sigma+1) a_{\infty}}{1-(\sigma-1) a_{\infty}}}B_2^{-2\frac{1-\sigma a_{\infty}}{1-(\sigma-1) a_{\infty}}}|t|^{\frac{1-(\sigma+1) a_{\infty}}{1-(\sigma-1) a_{\infty}}}\left(1+\mathcal{O}\left(\log|t||t|^{\frac{1-(\sigma+1) a_{\infty}}{1-(\sigma-1) a_{\infty}}}\right)\right),\\\label{bound:beta1-1}
    &|\Tilde{\beta}|=\mathcal{O}\left(|t|^{\frac{1-(\sigma+1)a_{\infty}}{1-(\sigma-1)a_{\infty}}}\right)\\
    \label{bound:lambda1-1}
    &\lambda(t)=B_2^{-\frac{1}{1+\frac{ a_{\infty}}{1-\sigma a_{\infty}}}} B_1^{\frac{1}{1+\frac{a_{\infty}}{1-\sigma a_{\infty}}}}|t|^{\frac{1}{1+\frac{a_{\infty}}{1-\sigma a_{\infty}}}}\left(1+\mathcal{O}\left(\log|t||t|^{\frac{1-(\sigma+1) a_{\infty}}{1-(\sigma-1) a_{\infty}}}\right)\right),\\\label{bound:alpha1-1}
    &\alpha(t)=b\lambda^{\frac{ a_{\infty}}{1-\sigma a_{\infty}}}=b_{\infty}B_2^{\frac{\frac{a_{\infty}}{1-\sigma a_{\infty}}}{1+\frac{a_{\infty}}{1-\sigma a_{\infty}}}}B_1^{\frac{\frac{a_{\infty}}{1-\sigma a_{\infty}}}{1+\frac{a_{\infty}}{1-\sigma a_{\infty}}}}|t|^{\frac{\frac{a_{\infty}}{1-\sigma a_{\infty}}}{1+\frac{a_{\infty}}{1-\sigma a_{\infty}}}}\left(1+\mathcal{O}\left(\log|t||t|^{\frac{1-(\sigma+1) a_{\infty}}{1-(\sigma-1) a_{\infty}}}\right)\right),\\\label{bound:gamma1-1}
    &\gamma(t)=(1+\beta_{\infty}^2)B_1^{-\frac{1-\frac{a_{\infty}}{1-\sigma a_{\infty}}}{1+\frac{a_{\infty}}{1-\sigma a_{\infty}}}}B_2^{2\frac{1}{1+\frac{a_{\infty}}{1-\sigma a_{\infty}}}}|t|^{-\frac{1-\frac{a_{\infty}}{1-\sigma a_{\infty}}}{1+\frac{a_{\infty}}{1-\sigma a_{\infty}}}}+\mathcal{O}(|\log t|),
\end{align}
where
\begin{align*}
    B_1=\frac{1-(\sigma-1) a_{\infty}}{1-(\sigma+1) a_{\infty}},~~B_2=\frac{ a_{\infty}}{2(1-(2\sigma+1) a_{\infty})\beta_{\infty}}b_{\infty}
\end{align*}
and a universal constant
\begin{align*}
   |b_{\infty}-1|\ll1.
\end{align*}
\end{lemma}
To prove Lemma \ref{lemma:initial-1}, we need to prove the following lemma.
\begin{lemma}\label{lemma:auxi-1}
There exists a universal constant $s_0\gg1$ such that the following holds. Let
\begin{align}\label{A:intial:a-1}
    \frac{1}{2}<b_0<1,~~\gamma_0\in\mathbb{R},~~a_0=\frac{1}{(1-(2\sigma+1) a_{\infty})s_0}.
\end{align}
Then the solution $(\lambda,a,\Tilde{\beta},\alpha,\gamma)$ to the dynamical system
\begin{align}\label{A:ode-1}
    \begin{cases}
    \frac{\lambda_s}{\lambda}=-a+\mathbf{A}_1(\mathcal{P}),~~\mathcal{P}=(a,\Tilde{\beta}),\\
    \frac{\alpha_s}{\lambda}=-2\beta,~~\beta=\beta_\infty+\Tilde{\beta},\\
    \beta_s=\mathbf{A}_2,\\
    \gamma_s=1+\beta^2,\\
    \frac{ds}{dt}=\frac{1}{\lambda^2},\\
    a=\frac{2\beta}{ a_{\infty}}\frac{\lambda}{\alpha},
    \end{cases}
    ~~\text{with}~~
    \begin{cases}
    \frac{\alpha(s_0)}{\lambda^{\frac{a_{\infty}}{1-\sigma a_{\infty}} }(s_0)}=b_0,\\
    a(s_0)=a_0,\\
    \Tilde{\beta}(s_0)=\frac{1}{s_0^2},\\
    \gamma(s_0)=\gamma_0,
    \end{cases}
\end{align}
is defined on $[s_0,+\infty)$. Moreover, there exists $b_\infty>0$ with
\begin{align}\label{bound:b-1}
    b_\infty=b_0+o_{s_{\infty}}(1)
\end{align}
such that the following asymptotics hold on $[s_0,+\infty)$:
\begin{align}\label{bound:a-1}
    &a(s)=\frac{1}{(1-(2\sigma+1) a_{\infty})s}+\mathcal{O}\left(\frac{|\log s|}{s^2}\right),~~|\Tilde{\beta}(s)|\lesssim\frac{1}{s},\\\label{bound:lambda-1}
    &\lambda(s)=\left(\frac{ a_{\infty}}{2(1-(2\sigma+1) a_{\infty})\beta_{\infty}s}b_{\infty}\right)^{\frac{1}{1-\frac{ a_{\infty}}{1-\sigma a_{\infty}}}}\left(1+\mathcal{O}\left(\frac{\log s}{s}\right)\right),\\\label{bound:alpha-1}
    &\alpha(t)=b\lambda^{\frac{ a_{\infty}}{1-\sigma a_{\infty}}}=b_{\infty}\left(\frac{ a_{\infty}}{2(1-(2\sigma+1)a_{\infty})\beta_{\infty}s}b_{\infty}\right)^{\frac{\frac{ a_{\infty}}{1-\sigma a_{\infty}}}{1-\frac{ a_{\infty}}{1-\sigma a_{\infty}}}}\left(1+\mathcal{O}\left(\frac{\log s}{s}\right)\right),\\\label{bound:gamma-1}
    &\gamma(s)=(1+\beta_{\infty}^2)s+\mathcal{O}(|\log s|).
\end{align}
\end{lemma}
\begin{proof}
{\bf Step 1:} Bootstrap bounds. From the Cauchy-Lipschitz theorem, we can obtain the local existence. To control the solution on large positive times, let us introduce the auxiliary function
\begin{align*}
    b=\frac{\alpha}{\lambda^{\frac{ a_{\infty}}{1-\sigma a_{\infty}}}},
\end{align*}
which from \eqref{A:ode-1} satisfies
\begin{align*}
    \frac{db}{ds}=&\frac{\alpha_s}{\lambda^{\frac{ a_{\infty}}{1-a_{\infty}\sigma}}}-\frac{ a_{\infty}\lambda_s\alpha}{(1-\sigma a_{\infty})\lambda^{\frac{ a_{\infty}}{1-\sigma a_{\infty}}+1}}=\frac{\alpha}{\lambda^{\frac{ a_{\infty}}{1-\sigma a_{\infty}}}}\frac{\alpha_s}{\alpha}-\frac{ a_{\infty}\lambda_s\alpha}{(1-\sigma a_{\infty})\lambda^{\frac{ a_{\infty}}{1-\sigma a_{\infty}}+1}}\notag\\
    =&-\frac{ a_{\infty}}{1-\sigma a_{\infty}}\frac{\alpha}{\lambda^{\frac{ a_{\infty}}{1-\sigma a_{\infty}}}}\left(\frac{\lambda_s}{\lambda}+(1-\sigma a_{\infty})a\right)=-\frac{ a_{\infty}}{1+\sigma a_{\infty}}b\left(\mathbf{A}_1-\sigma a_{\infty}a\right).
\end{align*}
Hence, system \eqref{A:ode-1} is equivalent to
\begin{align}\label{A:ode:2-1}
    \begin{cases}
    \frac{db}{ds}=-\frac{ a_{\infty}}{1-\sigma a_{\infty}}b\left(\mathbf{A}_1-\sigma a_{\infty}a\right),\\
    a_s+(1- a_{\infty})a^2=\frac{a}{\beta}\mathbf{A}_2+a\mathbf{A}_1,\\
    \Tilde{\beta}_s=\mathbf{A}_2,\\
    a=\frac{2\beta}{ a_{\infty}}\frac{\lambda}{\alpha},~~\beta=\beta_{\infty}+\Tilde{\beta},
    \end{cases}
    ~~\text{with}~~
    \begin{cases}
    a(s_0)=a_0,\\
    b(s_0)=b_0,\\
    \Tilde{\beta}(s_0)=\frac{1}{s_0^2}.
    \end{cases}
\end{align}
We assume the following priori bounds: for $s_0\leq s\leq s_1$,
\begin{align}\label{A:priori-1}
    |b(s)|\leq 1+2b_0,~~|\Tilde{\beta}(s)|\leq\frac{1}{s},~~\left|a(s)-\frac{1}{(1-(2\sigma+1) a_{\infty})s}\right|\leq\frac{(\log s)^2}{s^2}.
\end{align}

{\bf Step 2.} Closing the bootstrap. We calim that the bounds \eqref{A:priori-1} can be improved on $[s_0,s_1]$ provided $s_0$ has been chosen large enough. Indeed, let us close the $a$ bound. From \eqref{A:ode:2-1} and \eqref{A:priori-1}, we have
\begin{align*}
    \left|-\frac{d}{ds}\frac{1}{a}+(1- a_{\infty}-2\sigma a_{\infty})\right|=\left|\frac{a_s+(1- a_{\infty}-2\sigma a_{\infty})a^2}{a^2}\right|\\
    =\frac{\left|\frac{a}{\beta}\mathbf{A}_2+a\mathbf{A}_1-2\sigma a_{\infty} a^2\right|}{a^2}\lesssim\frac{a(a^2+|\Tilde{\beta}|^2)}{a^2}\lesssim\frac{1}{s},
\end{align*}
where in the Penultimate step we use the definition of $\mathbf{A}_1$ and $\mathbf{A}_2$ (see \eqref{def:a1a2-1}) and thus using the boundary condition on $a$ at $s_0$ and the initialization \eqref{A:intial:a-1}, we have
\begin{align*}
    \left|\frac{1}{a(s)}-(1- a_{\infty}-2\sigma a_{\infty})s\right|\lesssim\left|\frac{1}{a(s_0)}-(1- a_{\infty}-2\sigma a_{\infty})s_0\right|+\int_{s_0}^s\frac{1}{s}\lesssim \log s,
\end{align*}
where we assume $s\geq s_0\gg1$. Hence,
\begin{align}\label{A:a-1}
    \left|a(s)-\frac{1}{(1- a_{\infty}-2\sigma a_{\infty})s}\right|\lesssim\frac{\log s}{s^2}.
\end{align}

Next, we consider $\Tilde{\beta}$. Since
\begin{align*}
    \mathbf{A}_2=\sigma a_{\infty}\beta_0a+\mathcal{O}\left(a^2+|\Tilde{\beta}|^2\right),
\end{align*}
we obtain,
\begin{align*}
    \title{\beta}_s=\sigma a_{\infty}\beta_0a+\mathcal{O}\left(\frac{1}{s^2}\right).
\end{align*}
From \eqref{A:priori-1} and \eqref{A:a-1}, we have
\begin{align*}
    \left|\frac{d}{ds}\left(s^{\frac{2}{1- a_{\infty}-2\sigma a_{\infty}}}\Tilde{\beta}\right)\right|\lesssim s^{\frac{2}{1- a_{\infty}-2\sigma a_{\infty}}}\left(\frac{\log s}{s^2}|\Tilde{\beta}|+\frac{1}{s^2}\right)\lesssim s^{\frac{2}{1- a_{\infty}-2\sigma a_{\infty}}-2}.
\end{align*}
Using the boundary condition \eqref{A:ode:2-1} and $\frac{2}{1- a_{\infty}-2\sigma a_{\infty}}-2>0$, we get
\begin{align*}
    \left|s^{\frac{2}{1- a_{\infty}-2\sigma a_{\infty}}}\Tilde{\beta}(s)\right|\lesssim \left|s^{\frac{2}{1- a_{\infty}-2\sigma a_{\infty}}}\Tilde{\beta}(s_0)\right|+s^{\frac{2}{1- a_{\infty}-2\sigma a_{\infty}}-1}-s_0^{\frac{2}{1- a_{\infty}-2\sigma a_{\infty}}-1}\lesssim s^{\frac{2}{1- a_{\infty}-2\sigma a_{\infty}}-1},
\end{align*}
and thus
\begin{align}\label{A:beta-1}
    |\Tilde{\beta}(s)|\lesssim\frac{1}{s}.
\end{align}

Finally, we estimate $b$. From \eqref{A:ode:2-1}, \eqref{A:a-1} and \eqref{A:beta-1}, we have
\begin{align*}
    \left|\frac{d}{ds}b\right|\lesssim\frac{1+2b_0}{s^2}.
\end{align*}
Thus, by \eqref{A:ode:2-1}, we deduce
\begin{align}\label{A:b-1}
    |b(s)|\leq b_0+C\frac{1+2b_0}{s_0}\leq\frac{1}{2}+\frac{3}{2}b_0,
\end{align}
for $s_0$ large enough. The bounds \eqref{A:a-1}, \eqref{A:beta-1} and \eqref{A:b-1} improve \eqref{A:priori-1}, and thus from a standard continuity argument, the bounds \eqref{A:a-1}, \eqref{A:beta-1} and \eqref{A:b-1} hold on $[s_0,\infty)$ and the solution is global.

{\bf Step 3.} Conclusion. From \eqref{A:a-1} and \eqref{A:beta-1}, the bound \eqref{bound:a-1} is proved. Moreover, from \eqref{A:ode:2-1}, we have
\begin{align*}
    \int_{s_0}^\infty\left|\frac{d}{ds}b\right|\lesssim \int_{s_0}^\infty\frac{1}{s^2}=o(1)~~\text{as}~~s_0\to\infty,
\end{align*}
and hence there exists $b_{\infty}$ satisfying \eqref{bound:b-1} such that for any $s\geq s_0$,
\begin{align}\label{A:b2-1}
    |b(s)-b_\infty|\lesssim\frac{1}{s}.
\end{align}
From \eqref{A:ode:2-1}, \eqref{A:priori-1} and \eqref{A:b2-1}, we get
\begin{align*}
    \lambda(s)=\frac{ a_{\infty}a}{2\beta}\alpha=\frac{ a_{\infty}}{2(1-(2\sigma+1) a_{\infty})\beta_{\infty}s}b_{\infty}\lambda^{\frac{ a_{\infty}}{1-\sigma a_{\infty}}}\left(1+\mathcal{O}\left(\frac{\log s}{s}\right)\right).
\end{align*}
This means that
\begin{align*}
    \lambda(s)=\left(\frac{ a_{\infty}}{2(1-(2\sigma+1) a_{\infty})\beta_{\infty}s}b_{\infty}\right)^{\frac{1}{1-\frac{ a_{\infty}}{1-\sigma a_{\infty}}}}\left(1+\mathcal{O}\left(\frac{\log s}{s}\right)\right).
\end{align*}
By \eqref{A:b2-1}, we obtain
\begin{align*}
    \alpha(t)=b\lambda^{\frac{ a_{\infty}}{1-\sigma a_{\infty}}}=b_{\infty}\left(\frac{ a_{\infty}}{2(1-(2\sigma+1) a_{\infty})\beta_{\infty}s}b_{\infty}\right)^{\frac{\frac{ a_{\infty}}{1-\sigma a_{\infty}}}{1-\frac{ a_{\infty}}{1-\sigma a_{\infty}}}}\left(1+\mathcal{O}\left(\frac{\log s}{s}\right)\right).
\end{align*}
Finally, it remains to estimate $\gamma$. In view of \eqref{dynamical:sys-1} and \eqref{A:beta-1}, we have
\begin{align*}
    \frac{d\gamma}{ds}=1+\beta_{\infty}^2+\mathcal{O}\left(\frac{1}{s}\right),
\end{align*}
Integration between $s_0$ and $s$, we have
\begin{align*}
    \gamma(s)=(1+\beta_{\infty}^2)s+O(|\log s|).
\end{align*}

This complete the proof of Lemma \ref{lemma:auxi-1}.
\end{proof}

\begin{proof}[\bf Proof of Lemma \ref{lemma:initial-1}.]
Notice that from \eqref{bound:lambda-1}, we have
\begin{align*}
    \int_{s_0}^\infty\lambda^2<+\infty.
\end{align*}
Thus, since $\frac{ds}{dt}=\frac{1}{\lambda^2(t)}$, the time of existence of the dynamical system in time $t$ is finite, and we may choose the origin of time $t$ such that the final time is $0$. Then, for all $t_0\leq t<0$, we have
\begin{align*}
    -t=\int_s^{\infty}\lambda^2,
\end{align*}
which together with \eqref{bound:lambda-1} yields
\begin{align}\label{A:s-1}
    \frac{1}{s}
    =&\left(\frac{1-(\sigma-1) a_{\infty}}{1-(\sigma+1) a_{\infty}}\right)^{\frac{1-(\sigma+1) a_{\infty}}{1-(\sigma-1) a_{\infty}}}\left(\frac{2(1-(2\sigma+1) a_{\infty})\beta_{\infty}}{ a_{\infty}b_{\infty}}\right)^{2\frac{1-\sigma a_{\infty}}{1-(\sigma-1) a_{\infty}}}\notag\\
    &\times|t|^{\frac{1-(\sigma+1) a_{\infty}}{1-(\sigma-1) a_{\infty}}}\left(1+\mathcal{O}\left(\log|t||t|^{\frac{1+(\sigma-1) a_{\infty}}{1+(\sigma+1) a_{\infty}}}\right)\right).
\end{align}
Injecting \eqref{A:s-1} into \eqref{bound:a-1}, \eqref{bound:lambda-1}, \eqref{bound:alpha-1} and \eqref{bound:gamma-1}, we get the estimates \eqref{bound:a1-1}, \eqref{bound:beta1-1}, \eqref{bound:lambda1-1}, \eqref{bound:alpha1-1} and \eqref{bound:gamma1-1}.
This concludes the proof of Lemma \ref{lemma:initial-1}.

\end{proof}

\subsection{Geometrical decomposition}
In this subsection, from the standard argument, we show that there exists a unique decomposition which relies on the implicit function theorem, the mass subcritical nondegeneracy $(\Lambda Q,Q)\neq0$ and the modulation parameters is sufficiently small.

\begin{lemma}\label{lemma:decom:unique-1}
There exists a universal constant $\delta>0$ such that the following holds. Let $u$ be a radially symmetric function of the form
\begin{align*}
    u(t,x)= \alpha_0^{\frac{\sigma}{2}}(t)\frac{1}{\lambda_0(t)}Q_{(a_0,\Tilde{\beta}_0)}\left(t,\frac{r-\alpha_0}{\lambda_0}\right)e^{i\gamma_0}+\Tilde{u}_0(t,x)
\end{align*}
with
\begin{align*}
    \lambda_0, \alpha_0>0,~~\beta_0=\beta_\infty+\Tilde{\beta}_0,~~a_0=\frac{2\beta_0}{ a_{\infty}}\frac{\lambda_0}{ \alpha_0},
\end{align*}
the a priori bound
\begin{align}\label{lower:priori-1}
    \frac{\alpha_0}{\lambda_0^{ \frac{a_{\infty}}{1-\sigma a_{\infty}}}}\geq1,
\end{align}
and
\begin{align}\label{small:para-1}
    0<|a_0|+|\Tilde{\beta}_0|+\|\Tilde{u}_0\|_{L^2}<\delta.
\end{align}
Then there exists a unique decomposition
\begin{align*}
    u(t,x)=\alpha_1^{\frac{\sigma}{2}}(t)\frac{1}{\lambda_1(t)}Q_{(a_1,\Tilde{\beta}_1)}\left(t,\frac{r-\alpha_1}{\lambda_1}\right)e^{i\gamma_1}+\Tilde{u}_1(t,x)
\end{align*}
with
\begin{align*}
    \beta_1=\beta_{\infty}+\Tilde{\beta}_1,~~a_1=\frac{2\beta_1}{ a_{\infty}}\frac{\lambda_1}{\alpha_1},
\end{align*}
such that
\begin{align*}
    \Tilde{u}_1(x)=\alpha_1^{\frac{\sigma}{2}}(t)\frac{1}{\lambda_1(t)}\Tilde{\epsilon}_{(a_1,\Tilde{\beta}_1)}\left(\frac{r-\alpha_1}{\lambda_1}\right)e^{i\gamma_1}
\end{align*}
satisfies the orthogonality conditions
\begin{align*}
    (\Re\Tilde{\epsilon}_1,\xi_{a_1}yQ)=(\Re\Tilde{\epsilon}_1,\xi_{a_1}Q)=(\Im{\Tilde{\epsilon}_1},\xi_{a_1}\Lambda Q)=(\Im{\Tilde{\epsilon}}_1,\xi_{a_1}\partial_yQ)=0.
\end{align*}
Moreover, there holds the smallness
\begin{align*}
    \left|\frac{\lambda_1}{\lambda_0}-1\right|+\left|\frac{ \alpha_0-\alpha_1}{\lambda_0}\right|+|\Tilde{\beta}_0-\Tilde{\beta}_0|+|\gamma_0-\gamma_1|+\|\Tilde{u}_1\|_{L^2}\lesssim\delta.
\end{align*}
\end{lemma}
This is a standard consequence of the Implicit function theorem, For the convenience of readers', we will provide proof in Appendix \ref{appendix:unique-1}.

\subsection{Modulation equations and estimates}

Let $u(t,x)\in H^1$ be the radially solution to \eqref{equ:3-1} on a time interval $[t_0,t_1)$, $t_1<0$. From Lemma \ref{lemma:decom:unique-1},  we assume that $u(t)$ admits on $[t_0,t_1)$ a unique decomposition
\begin{align*}
    u(t,x)=\alpha^{\frac{\sigma}{2}}(t)\frac{1}{\lambda(t)}v\left(t,\frac{r-\alpha(t)}{\lambda(t)}\right)e^{i\gamma(t)},
\end{align*}
where
\begin{align}\label{def:a:ds-1}
    a(t)=\frac{2\beta}{ a_{\infty}}\frac{\lambda}{\alpha},~~\frac{ds}{dt}=\frac{1}{\lambda^2(t)},
\end{align}
and holds the unique decomposition
\begin{align}\label{decom:two-1}
    w(s,y)=v(s,y)e^{i\beta y}=\mathbf{Q}_{\mathcal{P}}e^{i\Tilde\beta y}+\Tilde{\epsilon}(t,y),~~\Tilde{\epsilon}=\Tilde{\epsilon}_1+i\Tilde{\epsilon}_2,~~\Tilde{\epsilon}=\epsilon e^{i\Tilde{\beta}}
\end{align}
with the orthogonality conditions
\begin{align}\label{decom:orth-1}
    (\Tilde{\epsilon}_1,\xi_a yQ)=(\Tilde{\epsilon}_1,\xi_a Q)=(\Tilde{\epsilon}_2,\xi_a\Lambda Q)=(\Tilde{\epsilon}_2,\xi_a Q^\prime)=0.
\end{align}
From \eqref{def:a:ds-1}, we have
\begin{align}\label{def:as:dcom-1}
    &a_s+(1- a_{\infty})a^2-\frac{a}{\beta}\mathbf{A}_2-b\mathbf{A}_1\notag\\
    =&\frac{a}{\beta}(\Tilde{\beta}_s-\mathbf{A}_2)+a\left(\frac{\lambda_s}{\lambda}+a-\mathbf{A}_1\right)-\frac{ a_{\infty}}{2\beta}a^2\left(\frac{\alpha_s}{\lambda}+2\beta\right).
\end{align}
The modulation equations are a consequence of the orthogonality conditions \eqref{decom:orth-1} and require the derivation of the equation for $\Tilde{\epsilon}$. Recall the equation \eqref{equ:w-1} satisfied by $w$:
\begin{align*}
     &i\partial_sw+w_{yy}-w+\frac{1}{\left|1+\frac{ a_{\infty}a}{2\beta}y\right|^\sigma}|w|^2w+\frac{ a_{\infty}a}{\beta}\frac{1}{1+\frac{ a_{\infty}a}{2\beta}y}(w_y-i\beta w)+a(i\Lambda w+\beta yw)-\frac{i\sigma a_{\infty}}{2}aw\notag\\
    =&-\beta_syw+\left(\frac{\lambda_s}{\lambda}+a\right)(i\Lambda w+\beta yw)+\left(\frac{\alpha_s}{\lambda}+2\beta\right)\left(iw_y+\beta w-\frac{i\sigma a_{\infty}a}{4\beta}w\right)+(\Tilde{\gamma}_s-\beta^2)w.
\end{align*}
We inject the decomposition \eqref{decom:two-1}, which we rewrite using \eqref{def:Qcut-1}:
\begin{align*}
    w=\xi_aQ_{\mathcal{P}}e^{-ia\frac{y^2}{4}}+\Tilde{\epsilon}.
\end{align*}
into \eqref{equ:w-1}, using the formula \eqref{def:as:dcom-1}, \eqref{def:Psi-1} and the fact that $P_{\mathcal{P}}=Q+\mathcal{O}(ae^{-c|y|})$, then we obtain the following system
\begin{align}\label{system:1-1}
    \partial_s\Tilde{\epsilon}_1-M_-\Tilde{\epsilon}_2=&-\frac{ a_{\infty}a}{2\beta}\frac{2}{1+\frac{ a_{\infty}a}{2\beta}y}(\partial_y\Tilde{\epsilon}_2-\beta\Tilde{\epsilon}_1)-\Tilde{\beta}_sy\Tilde{\epsilon}_2+\left(\frac{\lambda_s}{\lambda}+a-\mathbf{A}_1\right)\Lambda Q\notag\\
    &+\frac{\lambda_s}{\lambda}(\Lambda \Tilde{\epsilon}_1+\beta y\Tilde{\epsilon}_2)+\left(\frac{\alpha_s}{\lambda}+2\Tilde{\beta}\right)\left(\partial_yQ+\partial_y\Tilde{\epsilon}_1-\frac{\sigma a_{\infty}a}{4\beta}(Q+\Tilde{\epsilon}_1)\right)\notag\\
    &+\Gamma\Tilde{\epsilon}_2-\Im{R(\Tilde{\epsilon})}+\mathcal{O}(a|\Tilde{\epsilon}|+a^l+a\mathbf{Mod})e^{-c|y|},\\\label{system:2-1}
    \partial_s\Tilde{\epsilon}_2+M_+\Tilde{\epsilon}_1=&-\frac{ a_{\infty}a}{2\beta}\frac{2}{1+\frac{ a_{\infty}a}{2\beta}y}(-\partial_y\Tilde{\epsilon}_1-\beta\Tilde{\epsilon}_2)+(\beta_s-\mathbf{A}_2)yQ+\Tilde{\beta}_sy\Tilde{\epsilon}_1\notag\\
    &-\beta\left(\frac{\lambda_s}{\lambda}+a-\mathbf{A}_1\right)yQ+\frac{\lambda_s}{\lambda}(\Lambda\Tilde{\epsilon}_2-\beta y\Tilde{\epsilon}_1)+\left(\frac{\alpha_s}{\lambda}+2\Tilde{\beta}\right)\left(\partial_y\Tilde{\epsilon}_2-\frac{\sigma a_{\infty}a}{4\beta}\Tilde{\epsilon}_2\right)\notag\\
    &-\Gamma(Q+\Tilde{\epsilon}_1)+\Re{R(\Tilde{\epsilon})}+\mathcal{O}(a|\Tilde{\epsilon}|+a^l+a\mathbf{Mod})e^{-c|y|},
\end{align}
where $(M_+,M_-)$ are small deformation of the linearized operator $(L_+,L_-)$ close to $Q$:
\begin{equation*}
    \begin{aligned}
    M_+(\Tilde{\epsilon})=&-\partial_y^2\Tilde{\epsilon}_1+\Tilde{\epsilon}_1-3\frac{1}{\left|1+\frac{ a_{\infty}a}{2\beta}y\right|^{\sigma}}Q^2\Tilde{\epsilon}_1,\\
    M_-(\Tilde{\epsilon})=&-\partial_y^2\Tilde{\epsilon}_2+\Tilde{\epsilon}_2-\frac{1}{\left|1+\frac{ a_{\infty}a}{2\beta}y\right|^{\sigma}}Q^2\Tilde{\epsilon}_2,
    \end{aligned}
\end{equation*}
and
\begin{align}\label{def:Gamma-1}
    \Gamma=(\Tilde{\gamma}_s-\beta^2)+\beta\left(\frac{\alpha_s}{\lambda}+2\beta\right)
\end{align}
and the nonlinear term is given by
\begin{equation}\label{def:nonlinear-1}
    R(\Tilde{\epsilon})=\frac{1}{\left|1+\frac{ a_{\infty}a}{2\beta}y\right|^{\sigma}}\left[(Q+\Tilde{\epsilon})(\Tilde{\epsilon}_1^2+\Tilde{\epsilon}_2^2)+2Q\Tilde{\epsilon}_1\Tilde{\epsilon}\right].
\end{equation}
We the define
\begin{align*}
    \mathbf{Mod}(t)=\left|\frac{\alpha_s}{\lambda}+2\beta\right|+|\Tilde{\gamma}-\Tilde{\beta}^2|+\left|\frac{\lambda_s}{\lambda}-a-\mathbf{A}_1\right|+|\Tilde{\beta}_s-\mathbf{A}_2|.
\end{align*}

To prove the modulation estimate, we define the renormalized weight on the Lebesgue measure
\begin{align*}
    \mu=\left(1+\frac{\lambda(t)}{\alpha(t)}y\right)^2=\left(1+\frac{a_{\infty }a}{2\beta}y\right)^2.
\end{align*}
and the weighted Sobolev norms
\begin{align*}
    \|\epsilon\|_{L^2_{\mu}}^2=\int|\epsilon|^2\mu,~~\|\epsilon\|_{H^1_{\mu}}^2=\int|\partial_y\epsilon|^2\mu+\int|\epsilon|^2\mu.
\end{align*}

We now give the following modulation estimate.
\begin{lemma}
There holds the bounds
\begin{align}\label{mod:esti:1-1}
    &|\mathbf{Mod}(t)|\lesssim a\|\Tilde{\epsilon}\|_{H_\mu^1}+a^l,\\\label{mod:esti:2-1}
    &\left|a_s+(1- a_{\infty})a^2-\frac{a}{\beta}\mathbf{A}_2-a\mathbf{A}_1\right|\lesssim a^2\|\Tilde{\epsilon}\|_{H^1_\mu}+a^{l+1}.
\end{align}
\end{lemma}
\begin{proof}
Except for dealing with the new terms that comes from  the inhomogeneous factor, the proof of this lemma is similar to \cite{RS2011JAMS,MRS2014Duke,L2022}. For the reader convenience, we give the proof of this lemma.

We divide the proof into the following steps.

{\bf Step 1:} Estimate $\left|\frac{\alpha_s}{\lambda}+2\beta\right|$.

We multiply the equation of $\Tilde{\epsilon}_1$ \eqref{system:1-1} by
 $\xi_ayQ$ and integrate by parts. Using the orthogonality conditions \eqref{decom:orth-1}, identity $L_-(yQ)=-2Q^\prime$ and the following relation
\begin{align}\label{mod:iden:1-1}
    (\partial_yQ,\xi_ayQ)=-\frac{1}{2}\|Q\|_{L^2}^2+\mathcal{O}(e^{-\frac{c}{\sqrt{a}}}),
\end{align}
we obtain
\begin{align}\label{mod:e1-1}
    \left|\frac{\alpha_s}{\lambda}+2\beta\right|\lesssim& a\|\Tilde{\epsilon}\|_{L^2_\mu}+\mathbf{Mod}(a+\|\Tilde{\epsilon}\|_{L^2_\mu})+a^l+\int |y|^{c_l}|R(\Tilde{\epsilon})|\eta_be^{-|y|}\notag\\
    \lesssim&a\|\Tilde{\epsilon}\|_{L^2_\mu}+\mathbf{Mod}(a+\|\Tilde{\epsilon}\|_{L^2_\mu})+a^l+\|\Tilde{\epsilon}\|_{L^2_{\mu}}^2,
\end{align}
where we  used \eqref{def:nonlinear-1}, the decay estimate \eqref{decay:P-1}, $\Tilde{\epsilon}$ is small, the following bound
\begin{align}\label{GN-1}
    \|\Tilde{\epsilon}\|_{L^\infty(y\geq-\frac{\delta}{a})}\leq\|\Tilde{\epsilon}^\prime\|_{L^2(y\geq-\frac{\delta}{a})}^{\frac{1}{2}}\|\Tilde{\epsilon}\|_{L^2(y\geq-\frac{\delta}{a})}^{\frac{1}{2}}\lesssim\|\Tilde{\epsilon}\|_{H^1_{\mu}}
\end{align}
and
\begin{align*}
    \int |y|^{c_l}|R(\Tilde{\epsilon})|\xi_ae^{-|y|}\lesssim&\int|y|^{c_k}\xi_a^{2}e^{-2|y|}|\Tilde
    \epsilon|^2+\int |\Tilde\epsilon|^3\xi_b\\
    \lesssim&\|\epsilon\|_{L^2_\mu}^2+\|\Tilde\epsilon\|_{L^\infty}\|\Tilde\epsilon\|_{L^2_\mu}^2\lesssim\|\Tilde\epsilon\|_{L^2_\mu}^2.
\end{align*}
Here we also used the decay estimates \eqref{decay:P-1} and \eqref{GN-1}.

{\bf Step 2:} Estimate $\left|\frac{\lambda_s}{\lambda}+a-\mathbf{A}_1\right|$.

We multiply the equation of $\Tilde{\epsilon}_1$ \eqref{system:1-1} by $\xi_aQ$ and integrate by parts. Using the orthogonality condition \ref{decom:orth-1}, $L_-Q=0$ and the relation
\begin{align}\label{mod:iden:2-1}
    (\xi_a\Lambda Q,Q)=\frac{1}{2}\int Q^2+\mathcal{O}(e^{-\frac{c}{\sqrt{a}}}),
\end{align}
we have
\begin{align}\label{mod:e2-1}
    \left|\frac{\lambda_s}{\lambda}+a-\mathbf{A}_1\right|\lesssim a\|\Tilde{\epsilon}\|_{L^2_\mu}+\mathbf{Mod}(a+\|\Tilde{\epsilon}\|_{L^2_\mu})+a^l+\|\Tilde{\epsilon}\|_{L^2_\mu}^2.
\end{align}

{\bf Step 3:} Estimate $\Gamma$, where $\Gamma$ is defined by \eqref{def:Gamma-1}.

Multiplying the equation of $\Tilde{\epsilon}_2$ \eqref{system:2-1} by $\xi_a\Lambda Q$ and using the orthogonality condition \eqref{decom:orth-1}, $L_+\Lambda Q=-2Q$ and \eqref{mod:iden:2-1}, we deduce
\begin{align}\label{mod:e3-1}
    |\Gamma|\lesssim a\|\Tilde{\epsilon}\|_{L^2_\mu}+\mathbf{Mod}(a+\|\Tilde{\epsilon}\|_{L^2_\mu})+a^l+\|\Tilde{\epsilon}\|_{L^2_\mu}^2.
\end{align}

{\bf Step 4:} Estimate $|\Tilde{\beta}_s-\mathbf{A}_2|$.

We multiply the equation of $\Tilde{\epsilon}_2$ \eqref{system:2-1} by $\xi_a\partial_yQ$ and use the orthogonality condition \eqref{decom:orth-1}, $L_+(\partial_yQ)=0$ and \eqref{mod:iden:1-1} to obtain
\begin{align}\label{mod:e4-1}
    |\Tilde{\beta}_s-\mathbf{A}_2|\lesssim a\|\Tilde{\epsilon}\|_{L^2_\mu}+\mathbf{Mod}(a+\|\Tilde{\epsilon}\|_{L^2_\mu})+a^l+\|\Tilde{\epsilon}\|_{L^2_\mu}^2.
\end{align}
Combining the above estimates \eqref{mod:e1-1}, \eqref{mod:e2-1}, \eqref{mod:e3-1} and \eqref{mod:e4-1} yields \eqref{mod:esti:1-1}. Estimate \eqref{mod:esti:2-1} follows from \eqref{mod:esti:1-1} and \eqref{def:as:dcom-1}.
We now complete the proof of this lemma.
\end{proof}

\section{Refine energy estimate}
In this section, our aim is to derive a mixed energy/Morawetz type estimate which is very crucial for the blowup solution.

Let $u(t,x)$ be a solution to \eqref{equ:3-1} on $[t_0,0)$ and have the following decomposition
\begin{align*}
    u(t,x)=\Tilde{Q}(t,x)+\Tilde{u}(t,x),
\end{align*}
where $\Tilde{Q}(t,x)$ defined by \eqref{def:Qtilde-1} and
\[\Tilde{u}(t,x)=\alpha^{\frac{\sigma}{2}}(t)\frac{1}{\lambda(t)}\epsilon\left(s,\frac{r-\alpha(t)}{\lambda(t)}\right)e^{i\gamma(t)},~~\Tilde{\epsilon}(s,y)=\epsilon(s,y)e^{i\beta y}\]
which in view of \eqref{equ:Qtilde-1}, yields the equation for $\Tilde{u}$,
\begin{align}\label{equ:utilde-1}
    i\partial_t\Tilde{u}+\Delta\Tilde{u}+\frac{1}{|x|^{\sigma}}\left(|u|^2u-|\Tilde{Q}|^2\Tilde{Q}\right)=-\psi=-\alpha^\frac{\sigma}{2}(t)\frac{1}{\lambda^{3}(t)}\Psi\left(t,\frac{r-\alpha(t)}{\lambda(t)}\right)e^{i\gamma(t)},
\end{align}
with $\Psi$ defined by \eqref{def:Psi-1}. From the Lemma \ref{lemma:decom:unique-1}, the decomposition holds as long as
\[\frac{\alpha(t)}{\lambda^{\frac{a_{\infty}}{1-\sigma a_{\infty}}}}\gtrsim1\]
and
\[|a(t)|+|\Tilde{\beta}(t)|+\|\Tilde{\epsilon}(t)\|_{L^2_{\mu}}<\delta.\]
then we assume a priori bounds
\begin{align}\label{priori:1-1}
    \|\epsilon\|_{H^1_\mu}<\min(a,\lambda)\delta,~~0<a<\delta,~~|\Tilde{\beta}|\leq a,
\end{align}
and
\begin{align}\label{priori:2-1}
    \frac{b_0}{2}\leq\frac{\alpha(t)}{\lambda(t)^{\frac{a_{\infty}}{1-\sigma a_{\infty}}}}\leq2b_0,
\end{align}
where $b_0$ is defined in Lemma \ref{lemma:initial-1}.

Let $\phi:[-1,+\infty)\to\mathbb{R}$ be a time-independent smooth compactly supported cutoff function which satisfies
\begin{align}\label{def:cutoff:phi-1}
    \phi(z)=\begin{cases}
    0~~\text{for}~~-1\leq z\leq-\frac{1}{2}~~\text{and for}~~z\geq\frac{1}{2},\\
    1~~\text{for}~~z\in\left(-\frac{1}{2},\frac{1}{2}\right),
    \end{cases}
\end{align}
and
\begin{align*}
    \sup_{z\geq-1}|\phi(z)|<2.
\end{align*}
Let
\[F(u)=\frac{1}{4}|u|^4,~~f(u)=|u|^2u,~~F^\prime(u)\cdot h=\Re{f(u)\Bar{h}}.\]
Now we give the following energy/Virial estimate.
\begin{lemma}\label{Lemma:refine-1}
Let
\begin{align}\label{def:J-1}
    J(\Tilde{u})=&\frac{1}{2}\int|\nabla\Tilde{u}|^2+\frac{1+\beta^2}{2}\int\frac{|\Tilde{u}|^2}{\lambda^2}-\int \frac{1}{|x|^{\sigma}}[F(\Tilde{Q}+\Tilde{u})-F(\Tilde{Q})-F^\prime(\Tilde{Q})\cdot\Tilde{u}]\notag\\
    &+\frac{\beta}{\lambda}\Im{\int \phi\left(\frac{r}{\alpha(t)}-1\right)\partial_r\Tilde{u}\Tilde{u}},
\end{align}
and
\begin{align}\label{def:Ku-1}
    K(\Tilde{u})=&-\frac{1+\beta^2}{2}\Im\left(\frac{1}{|x|^{\sigma}}(f(u)-f(\Tilde{Q})),\Bar{\Tilde{u}}\right)-\frac{2\beta}{\lambda}\Re\left(\int\phi\left(\frac{r}{\alpha(t)}-1\right)\frac{1}{|x|^{\sigma}}(f(\Tilde{Q}+\Tilde{u})-f(\Tilde{Q}))\overline{\partial_r\Tilde{u}}\right)\notag\\
    &-\Re\left(\partial_t\Tilde{Q},\overline{\frac{1}{|x|^{\sigma}}(f(\Tilde{Q}+\Tilde{u})-f(\Tilde{Q})-f^\prime(\Tilde{Q})\cdot\Tilde{u}})\right).
\end{align}
Then the following holds:
\begin{align}\label{energy:virial-1}
    \frac{d}{dt}J(\Tilde{u})=K(\Tilde{u})+\mathcal{O}\left(\frac{a}{\lambda^4}\|\epsilon\|_{H^1_\mu}^2+\frac{a^l}{\lambda^4}\|\epsilon\|_{H^1_\mu}\right).
\end{align}
\end{lemma}
\begin{proof}
We divide the proof into the two steps.

{\bf Step 1:} Estimate the energy part. From \eqref{equ:utilde-1}, we have
\begin{align}\label{ref:energy:1-1}
    &\frac{d}{dt}\left\{\frac{1}{2}\int|\nabla\Tilde{u}|^2+\frac{1+\beta^2}{2}\int\frac{|\Tilde{u}|^2}{\lambda^2}-\int \frac{1}{|x|^{\sigma}}[F(\Tilde{Q}+\Tilde{u})-F(\Tilde{Q})-F^\prime(\Tilde{Q})\cdot\Tilde{u}]\right\}\notag\\
    =&-\Re\left(\partial_t\Tilde{u},\overline{\Delta\Tilde{u}-\frac{1+\beta^2}{\lambda^2}\Tilde{u}+\frac{1}{|x|^{\sigma}}(f(u)-f(\Tilde{Q}))}\right)-\frac{(1+\beta^2)\lambda_t}{\lambda^3}\int|\Tilde{u}|^2\notag\\
    &+\frac{\beta\beta_t}{\lambda^2}\int|\Tilde{u}|^2-\Re\left(\partial_t\Tilde{Q},\overline{\frac{1}{|x|^{\sigma}}(f(u)-f(\Tilde{Q})-f^\prime(\Tilde{Q})\cdot\Tilde{u})}\right)\notag\\
    =&\Im\left(\psi,\overline{\Delta\Tilde{u}-\frac{1+\beta^2}{\lambda^2}\Tilde{u}+\frac{1}{|x|^{\sigma}}(f(u)-f(\Tilde{Q}))}\right)-\frac{1+\beta^2}{\lambda^2}\Im\left(\frac{1}{|x|^{\sigma}}(f(u)-f(\Tilde{Q})),\Bar{\Tilde{u}}\right)\notag\\
    &-\frac{(1+\beta^2)\lambda_t}{\lambda^3}\int|\Tilde{u}|^2+\frac{\beta\beta_t}{\lambda^2}\int|\Tilde{u}|^2\notag\\
    &-\Re\left(\partial_t\Tilde{Q},\overline{\frac{1}{|x|^{\sigma}}(f(u)-f(\Tilde{Q})-f^\prime(\Tilde{Q})\cdot\Tilde{u})}\right).
\end{align}
From \eqref{mod:esti:1-1} and the identity $\frac{\lambda_t}{\lambda^3}=\frac{\lambda_s}{\lambda^5}=\frac{a}{\lambda^4}-\frac{\mathbf{A}_1}{\lambda^4}-\frac{1}{\lambda^4}\left(\frac{\lambda_s}{\lambda}+a-\mathbf{A}_1\right)$, we have
\begin{align}\label{ref:energy:2-1}
    -\frac{\lambda_t}{\lambda^3}\int|\Tilde{u}|^2=&\frac{a}{\lambda^4}\int|\Tilde{u}|^2-\frac{\mathbf{A}_1}{\lambda^4}\int|\Tilde{u}|^2-\frac{1}{\lambda^4}\left(\frac{\lambda_s}{\lambda}+a-\mathbf{A}_1\right)\int|\Tilde{u}|^2\notag\\
    =&\frac{1}{\lambda^4}\mathcal{O}(a\|\epsilon\|_{L^2_\mu}^2),
\end{align}
where we used the priori estimate \eqref{priori:1-1}, \eqref{priori:2-1} and Lemma \ref{lemma:initial-1}. Also from \eqref{mod:esti:1-1}, we have
\begin{align}\label{ref:energy:3-1}
    \frac{\beta\beta_t}{\lambda^2}\int|\Tilde{u}|^2=\frac{\beta\mathbf{A}_2}{\lambda^4}\int|\Tilde{u}|^2+\frac{\beta(\beta_s-\mathbf{A}_2)}{\lambda^4}\int|\Tilde{u}|^2
    =\frac{1}{\lambda^4}\mathcal{O}(a\|\epsilon\|_{L^2_\mu}^2).
\end{align}
From \eqref{def:Psi-1}, \eqref{mod:esti:1-1} and \eqref{mod:esti:2-1}, we have
\begin{align}\label{ref:energypsi-1}
    |\psi|\lesssim& \xi_a(a^l+\mathbf{Mod})(1+|y|^{c_l})e^{-|y|}+\frac{e^{-|y|}}{a^{c_l}}\mathbf{1}_{y\sim\frac{1}{\sqrt{a}}}\notag\\
    \lesssim&\xi_a(a^l+a\|\epsilon\|_{H^1_\mu})(1+|y|^{c_l})e^{-|y|}+\frac{e^{-|y|}}{a^{c_l}}\mathbf{1}_{y\sim\frac{1}{\sqrt{a}}}.
\end{align}
Now we estimate the first term in the right-hand side of \eqref{ref:energy:1-1}, using the fact that $f^\prime(w)\cdot u=2|w|^2u+w^2\Bar{u}$ and the priori estimate \eqref{priori:1-1} and \eqref{priori:2-1} and \eqref{ref:energypsi-1}, we have
\begin{align}\label{ref:energy:4-1}
    &\left|\Im\left(\psi,\overline{\Delta\Tilde{u}-\frac{1+\beta^2}{\lambda^2}\Tilde{u}+\frac{1}{|x|^{\sigma}}(f(u)-f(\Tilde{Q}))}\right)\right|\notag\\
    =&\left|\Im\left(\psi,\overline{\Delta\Tilde{u}-\frac{1+\beta^2}{\lambda^2}\Tilde{u}+\frac{1}{|x|^{\sigma}}(f(u)-f(\Tilde{Q})-f^{\prime}(\Tilde{Q})\cdot\Tilde{u}+f^{\prime}(\Tilde{Q})\cdot\Tilde{u}))}\right)\right|\notag\\
    \lesssim&\left|\Im\int\left(\Delta\psi-\frac{1+\beta^2}{\lambda^2}\psi+2\frac{1}{|x|^{\sigma}}|\Tilde{Q}|^2\psi+\frac{1}{|x|^{\sigma}}\Tilde{Q}^2\Bar{\psi}\right)\Bar{\Tilde{u}}\right|+\left|\Im\left(\psi,\overline{\frac{1}{|x|^{\sigma}}(f(u)-f(\Tilde{Q})-f^{\prime}(\Tilde{Q})\cdot\Tilde{u})}\right)\right|\notag\\
    \lesssim&\frac{\left(a^l+a\|\epsilon\|_{H^1_\mu}\right)}{\lambda^4}\|\epsilon\|_{H^1_\mu}+\frac{1}{\lambda^4}\int \left(\xi_a(a^l+\mathbf{Mod})(1+|y|^{c_l})e^{-|y|}+\frac{e^{-|y|}}{a^{c_a}}\mathbf{1}_{y\sim\frac{1}{\sqrt{a}}}\right)\frac{1}{|x|^{\sigma}}\left||\Tilde{Q}||\Tilde{u}|^2+|\Tilde{u}|^3\right|\notag\\
    \lesssim&\frac{\left(a^l+a\|\epsilon\|_{H^1_\mu}\right)}{\lambda^4}\|\epsilon\|_{H^1_\mu}
\end{align}
where in the penultimate and last step, we used the Sobolev embedding and Lemma \ref{lemma:initial-1}.
Injecting \eqref{ref:energy:2-1}, \eqref{ref:energy:3-1} and \eqref{ref:energy:4-1}  into \eqref{ref:energy:1-1} yields
\begin{align}\label{energy:part-1}
    &\frac{d}{dt}\left\{\frac{1}{2}\int|\nabla\Tilde{u}|^2+\frac{1+\Tilde{\beta}^2}{2}\int\frac{|\Tilde{u}|^2}{\lambda^2}-\int \frac{1}{|x|^{\sigma}}[F(\Tilde{Q}+\Tilde{u})-F(\Tilde{Q})-F^\prime(\Tilde{Q})\cdot\Tilde{u}]\right\}\notag\\
    =&-\frac{1+\beta^2}{\lambda^2}\Im\left(\frac{1}{|x|^{\sigma}}(f(u)-f(\Tilde{Q})),\Bar{\Tilde{u}}\right)-\Re\left(\partial_t\Tilde{Q},\overline{\frac{1}{|x|^{\sigma}}(f(u)-f(\Tilde{Q})-f^\prime(\Tilde{Q})\cdot\Tilde{u})}\right)\notag\\
    &+\frac{1}{\lambda^4}\mathcal{O}(a^l\|\epsilon\|_{H^1_\mu}+a\|\epsilon\|_{H^1_\mu}^2).
\end{align}

{\bf Step 2:} Estimate the Virial part. We now estimate the other part in \eqref{def:J-1}. Using \eqref{def:a:ds-1}, we have the following relation
\begin{align*}
    \frac{d}{dt}\frac{r}{\alpha(t)}=-\frac{\alpha_tr}{\alpha^2(t)}=-\frac{\alpha_s}{\lambda}\frac{r}{\lambda \alpha^2(t)}=\frac{ a_{\infty}a(t)}{\lambda^2}\frac{r}{\alpha(t)}-\frac{ a_{\infty}a}{2\beta\lambda^2(t)}\frac{r}{\alpha(t)}\left(\frac{\alpha_s}{\lambda}+2\beta\right).
\end{align*}
From the above relation, by the direct computation, we can deduce the following
\begin{align}\label{virial:1-1}
    &\frac{d}{dt}\left\{\frac{\beta}{\lambda}\Im\int\phi\left(\frac{r}{\alpha(t)}-1\right)\partial_t\Tilde{u}\Bar{\Tilde{u}}\right\}\notag\\
    =&\frac{ a_{\infty}\beta a}{\lambda^3}\Im\int\frac{r}{\alpha(t)}\phi^\prime\left(\frac{r}{\alpha(t)}-1\right)\partial_r\Tilde{u}\Bar{\Tilde{u}}\notag\\
    &-\frac{ a_{\infty}a}{2\lambda^3}\left(\frac{\alpha_s}{\lambda}+2\beta\right)\Im\int\frac{r}{\alpha(t)}\phi^\prime\left(\frac{r}{\alpha(t)}-1\right)\partial_r\Tilde{u}\Bar{\Tilde{u}}\notag\\
    &+\frac{\mathbf{A}_2}{\lambda^3}\Im\int\phi\left(\frac{r}{\alpha(t)}-1\right)\partial_r\Tilde{u}\Bar{\Tilde{u}}+\frac{\beta_s-\mathbf{A}_2}{\lambda^3}\Im\int\phi\left(\frac{r}{\alpha(t)}-1\right)\partial_r\Tilde{u}\Bar{\Tilde{u}}\notag\\
    &+\frac{\beta(a-\mathbf{A}_1)}{\lambda^3}\Im\int\phi\left(\frac{r}{\alpha(t)}-1\right)\partial_r\Tilde{u}\Bar{\Tilde{u}}\notag\\
    &-\frac{\beta}{\lambda^3}\left(\frac{\lambda_s}{\lambda}+a-\mathbf{A}_1\right)\Im\int\phi\left(\frac{r}{\alpha(t)}-1\right)\partial_r\Tilde{u}\Bar{\Tilde{u}}\notag\\
    &+\frac{\beta}{\lambda}\Re\int i\partial_t\Tilde{u}\left(\overline{\partial_r\left(\phi\left(\frac{r}{\alpha(t)}-1\right)\right)\Tilde{u}+2\phi\left(\frac{r}{\alpha(t)}-1)\right)\partial_r\Tilde{u}}\right)\notag\\
    =&\frac{\beta}{\lambda}\Re\int i\partial_t\Tilde{u}\left(\overline{\partial_r\left(\phi\left(\frac{r}{\alpha(t)}-1\right)\right)\Tilde{u}+2\phi\left(\frac{r}{\alpha(t)}-1)\right)\partial_r\Tilde{u}}\right)+\mathcal{O}\left(\frac{a}{\lambda^4}\|\epsilon\|_{H^1_\mu}^2\right),
\end{align}
where we used \eqref{mod:esti:1-1}, Lemma \ref{lemma:initial-1} and the priori estimates \eqref{priori:1-1} and \eqref{priori:2-1}, and the fact that
\begin{align}\label{virial:phi-1}
    \frac{1}{r}\sim\frac{1}{\alpha(t)}~~\text{on the support of}~\phi\left(\frac{\cdot}{\alpha(t)}-1\right).
\end{align}
Now from \eqref{equ:utilde-1}, \eqref{virial:phi-1} and the integration by parts, the first term in the right-hand side of \eqref{virial:1-1} can be write as follows:
\begin{align}\label{virial:2-1}
    &\frac{\beta}{\lambda}\Re\int i\partial_t\Tilde{u}\left(\overline{\frac{1}{\alpha(t)}\phi^\prime\left(\frac{r}{\alpha(t)}-1\right)\Tilde{u}+2\phi\left(\frac{r}{\alpha(t)}-1)\right)\partial_r\Tilde{u}}\right)\notag\\
    =&\frac{ a_{\infty}a}{\lambda^2}\int\phi^\prime\left(\frac{r}{\alpha(t)}-1\right)|\partial_r\Tilde{u}|^2-\frac{ a_{\infty}^2a^2}{8\beta\lambda^3}\int \Delta\left(\partial_r\left(\phi\left(\frac{r}{\alpha(t)}-1\right)\right)\right)|\Tilde{u}|^2\notag\\
    &-\frac{2\beta}{\lambda}\Re\int \phi\left(\frac{r}{\alpha(t)}-1\right)\frac{1}{|x|^{\sigma}}(f(\Tilde{Q}+\Tilde{u})-f(\Tilde{Q}))\overline{\partial_r\Tilde{u}}\notag\\
    &-\frac{\beta}{\lambda}\Re\int \partial_r\left(\phi\left(\frac{r}{\alpha(t)}-1\right)\right)\frac{1}{|x|^{\sigma}}(f(\Tilde{Q}+\Tilde{u})-f(\Tilde{Q}))\overline{\Tilde{u}}\notag\\
    &-\frac{2\beta}{\lambda}\Re{\int \phi\left(\frac{r}{\alpha(t)}-1\right)\psi\partial_r\Bar{\Tilde{u}}}-\frac{\beta}{\lambda}\Re\int\partial_r\left(\phi\left(\frac{r}{\alpha(t)}-1\right)\right)\psi\Bar{\Tilde{u}}\notag\\
    =&-\frac{2\beta}{\lambda}\Re\int \phi\left(\frac{r}{\alpha(t)}-1\right)\frac{1}{|x|^{\sigma}}(f(\Tilde{Q}+\Tilde{u})-f(\Tilde{Q}))\overline{\partial_r\Tilde{u}}\notag\\
    &-\frac{\beta}{\lambda}\Re\int \partial_r\left(\phi\left(\frac{r}{\alpha(t)}-1\right)\right)\frac{1}{|x|^{\sigma}}(f(\Tilde{Q}+\Tilde{u})-f(\Tilde{Q}))\overline{\Tilde{u}}\notag\\
    &-\frac{2\beta}{\lambda}\Re{\int \phi\left(\frac{r}{\alpha(t)}-1\right)\psi\partial_r\Bar{\Tilde{u}}}-\frac{\beta}{\lambda}\Re\int\partial_r\left(\phi\left(\frac{r}{\alpha(t)}-1\right)\right)\psi\Bar{\Tilde{u}}+\mathcal{O}\left(\frac{a}{\lambda^4}\|\epsilon\|_{H^1_\mu}^2\right).
\end{align}
Now we estimate the second term in the right hand side of \eqref{virial:2-1}
\begin{align*}
    &\left|-\frac{\beta}{\lambda}\Re\int \partial_r\left(\phi\left(\frac{r}{\alpha(t)}-1\right)\right)\frac{1}{|x|^{\sigma}}(f(\Tilde{Q}+\Tilde{u})-f(\Tilde{Q}))\overline{\Tilde{u}}\right|\notag\\
    =&\left|-\frac{\beta}{\lambda}\Re\int \partial_r\left(\phi\left(\frac{r}{\alpha(t)}-1\right)\right)\frac{1}{|x|^{\sigma}}(f(\Tilde{Q}+\Tilde{u})-f(\Tilde{Q})-f^{\prime}(\Tilde{Q})\cdot\Tilde{u}+f^{\prime}(\Tilde{Q})\cdot\Tilde{u})\overline{\Tilde{u}}\right|\notag\\
    \lesssim&\left|-\frac{\beta}{\lambda}\Re\int \partial_r\left(\phi\left(\frac{r}{\alpha(t)}-1\right)\right)\frac{1}{|x|^{\sigma}}(f(\Tilde{Q}+\Tilde{u})-f(\Tilde{Q})-f^{\prime}(\Tilde{Q})\cdot\Tilde{u})\overline{\Tilde{u}}\right|\notag\\
    &+\left|-\frac{\beta}{\lambda}\Re\int \partial_r\left(\phi\left(\frac{r}{\alpha(t)}-1\right)\right)\frac{1}{|x|^{\sigma}}(f^{\prime}(\Tilde{Q})\cdot\Tilde{u})\overline{\Tilde{u}}\right|\notag\\
    \lesssim&\frac{a}{\lambda^2}\int\frac{1}{|x|^{\sigma}}\left( |\Tilde{u}|^4+|\Tilde{u}|^3|\Tilde{Q}|\right)+\frac{a}{\lambda^2}\int
\frac{1}{|x|^{\sigma}}|\Tilde{Q}|^2|\Tilde{u}|^2
\end{align*}
where we used the \eqref{virial:phi-1}. Notice that $\Tilde{Q}$ is localized in the region $r\geq\frac{\alpha(t)}{2}$ due to the cut-off function $\xi_a$ in its definition. Now, the region $r\geq\frac{\alpha(t)}{2}$ corresponds to $y\geq-\frac{\alpha(t)}{2\lambda(t)}$ and thus
\begin{align*}
    \mu\gtrsim1~~\text{for}~~r\geq\frac{\alpha(t)}{2}.
\end{align*}
By the Sobolev embedding and priori bound \eqref{priori:1-1} and Lemma \ref{lemma:initial-1}, we have
\begin{align*}
    \int\frac{1}{|x|^{\sigma}}|\Tilde{Q}||\Tilde{u}|^3=&\frac{1}{\lambda^3}\int_{y\geq-\frac{\alpha}{2\lambda}}|\alpha|^{2\sigma+2}|\lambda y+\alpha|^{-\sigma}|\epsilon|^3|Q_{\mathcal{P}}|\mu\\
    \lesssim&\frac{1}{\lambda^3}\int_{y\geq-\frac{\alpha}{2\lambda}}|\alpha|^{\sigma+2}\left|\frac{\lambda}{\alpha} y+1\right|^{-\sigma}|\epsilon|^3\mu\\
    \lesssim&\frac{1}{\lambda^2}\int_{y\geq-\frac{\alpha}{2\lambda}}\left|\frac{\lambda}{\alpha} y+1\right|^{-\sigma}|\epsilon|^3\mu\\
    \lesssim&\frac{1}{\lambda^2}\|\epsilon\|_{L^{\infty}(y\geq-\frac{\alpha(t)}{2\lambda(t)})}\|\epsilon\|_{L^2_{\mu}}^2\lesssim\frac{1}{\lambda^2}\|\epsilon\|_{H^1_{\mu}}\|\epsilon\|_{L^2_{\mu}}^2\\
    \lesssim&\frac{\delta\|\epsilon\|_{H^1_{\mu}}^2}{\lambda^2}.
\end{align*}
Similarly, we can obtain
\begin{align*}
    \int\frac{1}{|x|^{\sigma}}|\Tilde{Q}|^2|\Tilde{u}|^2\lesssim\frac{\|\epsilon\|_{H^1_{\mu}}^2}{\lambda^2}.
\end{align*}
By the Gagliardo-Nirebergy inequality, we have
\begin{align*}
    \int\frac{1}{|x|^{\sigma}}|\Tilde{u}|^4\lesssim\|\nabla\Tilde{u}\|_{L^2}^{3+\sigma}\|\Tilde{u}\|_{L^2}^{1-\sigma}=\frac{1}{\lambda^2}\|\nabla\epsilon\|_{L^2_{\mu}}^{3+\sigma}\|\epsilon\|_{L^2_{\mu}}^{1-\sigma}\lesssim\frac{\delta^2\|\epsilon\|_{H^1_{\mu}}^2}{\lambda^2}.
\end{align*}
From above estimates, we can obtain
\begin{align}\label{virial:3-1}
    \left|-\frac{\beta}{\lambda}\Re\int \partial_r\left(\phi\left(\frac{r}{\alpha(t)}-1\right)\right)\frac{1}{|x|^{\sigma}}(f(\Tilde{Q}+\Tilde{u})-f(\Tilde{Q}))\overline{\Tilde{u}}\right|
    \lesssim\mathcal{O}\left(\frac{a}{\lambda^4}\|\epsilon\|_{H^1_\mu}^2\right).
\end{align}
We estimate the terms that contain $\psi$ in \eqref{virial:2-1}. Using \eqref{ref:energypsi-1}, we obtain
\begin{align}\label{virial:psi-1}
    &\left|-\frac{2\beta}{\lambda}\Re{\int \phi\left(\frac{r}{\alpha(t)}-1\right)\psi\partial_r\Bar{\Tilde{u}}}-\frac{\beta}{\lambda}\Re\int\partial_r\left(\phi\left(\frac{r}{\alpha(t)}-1\right)\right)\psi\Bar{\Tilde{u}}\right|\notag\\
    \lesssim&\frac{a^l+a\|\epsilon\|_{H^1_\mu}}{\lambda^4}\|\epsilon\|_{H^1_\mu}.
\end{align}
Injecting \eqref{virial:2-1}, \eqref{virial:3-1} and \eqref{virial:psi-1} into \eqref{virial:1-1}, we can obtain
\begin{align}\label{virial:part-1}
     &\frac{d}{dt}\left\{\frac{\Tilde{\beta}}{\lambda}\Im\int\phi\left(\frac{r}{\alpha(t)}-1\right)\partial_t\Tilde{u}\Bar{\Tilde{u}}\right\}\notag\\
     =&-\frac{2\Tilde{\beta}}{\lambda}\Re\int \phi\left(\frac{r}{\alpha(t)}-1\right)\frac{1}{|x|^{\sigma}}(f(\Tilde{Q}+\Tilde{u})-f(\Tilde{Q}))\overline{\partial_r\Tilde{u}}+\mathcal{O}\left(\frac{a}{\lambda^4}\|\epsilon\|_{H^1_\mu}^2+\frac{a^l}{\lambda^4}\|\epsilon\|_{H^1_\mu}\right).
\end{align}
Now combining the \eqref{energy:part-1} and \eqref{virial:part-1}, we can obtain \eqref{energy:virial-1}. This conclude the proof of Lemma \ref{Lemma:refine-1}.
\end{proof}

Next Lemma we will give the following upper bound and lower bound of \eqref{def:J-1}.
\begin{lemma}
Let $J$ defined by \eqref{def:J-1}. Then
\begin{align}\label{Bound:J-1}
    c_1\left(\|\nabla\Tilde{u}\|_{L^2}^2+\frac{1}{\lambda^2}\|\Tilde{u}\|_{L^2}\right)\geq J(\Tilde{u})\geq c_2\left(\|\nabla\Tilde{u}\|_{L^2}^2+\frac{1}{\lambda^2}\|\Tilde{u}\|_{L^2}\right)
\end{align}
for some constants $c_1,c_2>0$.
\end{lemma}
\begin{proof}
To prove the lower bound. Notice that $J(\Tilde{u})$ can be write as
\begin{align}\label{coercivity:1-1}
    J(\Tilde{u})=\frac{1}{2\lambda^2}&\Bigg(\int|\partial_y\epsilon|^2\mu+2\beta\Im\int\phi(z)\partial_y\epsilon\Bar{\epsilon}\mu+(1+\Tilde{\beta}^2)\int|\epsilon|^2\mu\notag\\
    &-2\int \left|\frac{\lambda}{\alpha} y+1\right|^{-\sigma}\left(F(\mathbf{Q}_{\mathcal{P}}+\epsilon)-F(\mathbf{Q}_{\mathcal{P}})-F^\prime(\mathbf{Q}_{\mathcal{P}})\cdot\epsilon\right)\mu\Bigg),
\end{align}
where
\[z=\frac{r}{\alpha}-1=\frac{ a_{\infty}a}{2\beta}y,~~\mu=(1+z)^2\]
and we used Lemma \ref{lemma:initial-1}.

We now estimate the nonlinear term in \eqref{coercivity:1-1}.
\begin{align}\label{coer:nonlinear-1}
    &\int \left|\frac{\lambda}{\alpha} y+1\right|^{-\sigma}\left(F(\mathbf{Q}_{\mathcal{P}}+\epsilon)-F(\mathbf{Q}_{\mathcal{P}})-F^\prime(\mathbf{Q}_{\mathcal{P}})\cdot\epsilon\right)\mu)\notag\\
    =&\int \left|\frac{\lambda}{\alpha} y+1\right|^{-\sigma}\left[\frac{1}{2}|\mathbf{Q}_{\mathcal{P}}|^2|\epsilon|^2+(\Re(\mathbf{Q}_{\mathcal{P}}\bar{\epsilon}))^2+\frac{1}{4}|\epsilon|^4+|\epsilon|^2\Re(\mathbf{Q}_{\mathcal{P}}\Bar{\epsilon})\right]\mu\notag\\
    =&\frac{1}{2}\int \left|\frac{\lambda}{\alpha} y+1\right|^{-\sigma}\xi_aQ^2(3\Tilde{\epsilon}_1^2+\Tilde{\epsilon}_2^2)+\mathcal{O}(\delta^{C}\|\epsilon\|_{H^1_\mu}^2),
\end{align}
where we used the fact that
\[\mathbf{Q}_{\mathcal{P}}=\xi_aQe^{-i\beta y}+\mathcal{O}(ae^{-c|y|}),~~\text{and}~~\Tilde{\epsilon}=\epsilon e^{i\beta y}.\]
Combining the above estimate \eqref{coer:nonlinear-1}, priori estimate \eqref{priori:1-1} and the definition of $\Tilde{\beta}$, we can get
\begin{align}\label{coer:2-1}
    J(\Tilde{u})=\frac{1}{2\lambda^2}&\Bigg\{\int|\partial_y\epsilon|^2\mu+2\beta_{\infty}\Im\int\phi(z)\partial_y\epsilon\Bar{\epsilon}\mu+(1+\beta_{\infty}^2)\int|\epsilon|^2\mu\notag\\
    &-\int \left|\frac{\lambda}{\alpha} y+1\right|^{-\sigma}\xi_aQ^2(3\Tilde{\epsilon}_1^2+\Tilde{\epsilon}_2^2)+\mathcal{O}\left(\delta^{C_k}\|\epsilon\|_{H^1_\mu}^2\right)\Bigg\}\notag\\
    =\frac{1}{2\lambda^2}&\Bigg\{I_1+I_2-\int \left|\frac{\lambda}{\alpha}y+1\right|^{-\sigma}\xi_aQ^2(3\Tilde{\epsilon}_1^2+\Tilde{\epsilon}_2^2)+\mathcal{O}(\delta^{C_k}\|\epsilon\|_{H^1_\mu}^2)\Bigg\},
\end{align}
where
\begin{align*}
    I_1=\int_{|y|\geq\frac{1}{\sqrt{a}}}|\partial_y\epsilon|^2\mu+2\beta_{\infty}\Im\int_{|y|\geq\frac{1}{\sqrt{a}}}\phi(z)\partial_y\epsilon\Bar{\epsilon}\mu+(1+\beta_{\infty}^2)\int_{|y|\geq\frac{1}{\sqrt{a}}}|\epsilon|^2\mu,\\
    I_2=\int_{|y|\leq\frac{1}{\sqrt{a}}}|\partial_y\epsilon|^2\mu+2\beta_{\infty}\Im\int_{|y|\leq\frac{1}{\sqrt{a}}}\phi(z)\partial_y\epsilon\Bar{\epsilon}\mu+(1+\beta_{\infty}^2)\int_{|y|\leq\frac{1}{\sqrt{a}}}|\epsilon|^2\mu.
\end{align*}
By the definition of $\phi$ (see \eqref{def:cutoff:phi-1}) and the basic property of the quadratic form, we can obtain
\begin{align}\label{coer:I1-1}
    I_1\gtrsim \int_{_{|y|\geq\frac{1}{\sqrt{b}}}}(|\partial_y\epsilon|^2+|\epsilon|^2)\mu,
\end{align}
since $\beta_{\infty}^2\psi^2(z)-(1+\beta_{\infty}^2)^2<0$.

For the case $|y|\leq\frac{1}{\sqrt{a}}$, from \eqref{def:cutoff:phi-1}, then $|\phi(z)-1|\lesssim|z|\lesssim\sqrt{a}$, hence
\begin{align}\label{coer:I2-1}
    I_2=\int_{|y|\leq\frac{1}{\sqrt{a}}}(|\partial_y\epsilon|^2+|\epsilon|^2)\mu+\mathcal{O}(\sqrt{a}\|\epsilon\|_{H^1_\mu}^2).
\end{align}
Combining the above estimates \eqref{coer:I1-1}, \eqref{coer:I2-1} and \eqref{coer:2-1}, we have
\begin{align}\label{coer:3-1}
    2\lambda^2J(\Tilde{u})=&\int_{|y|\leq\frac{1}{\sqrt{a}}}|\partial_y\Tilde{\epsilon}|^2+|\Tilde{\epsilon}|^2-\int \left|\frac{\lambda}{\alpha}y+1\right|^{-\sigma}\xi_aQ^2(3\Tilde{\epsilon}_1^2+\Tilde{\epsilon}_2^2)\notag\\
    &+\int_{|y|\geq\frac{1}{\sqrt{a}}}(|\partial_y\epsilon|^2+|\epsilon|^2)\mu+\mathcal{O}(\delta^{C_k}\|\epsilon\|_{H^1_\mu}^2).
\end{align}

On the other hand, by using the following property of the linearized operator $L=(L_+,L_-)$ (see \cite{CGNT2007SIAM})
\begin{align*}
    (L_+\epsilon_1,\epsilon_1)+(L_-\epsilon_2,\epsilon_2)\geq C_0\|\epsilon\|_{H^1}^2-\frac{1}{C_0}\Big((\epsilon_1,Q)^2+(\epsilon,yQ)^2+(\epsilon_2,\Lambda Q)^2\Big),
\end{align*}
and the orthogonality conditions \eqref{decom:orth-1}, we can obtain
\begin{align}\label{coer:4-1}
    &\int_{|y|\leq\frac{1}{\sqrt{a}}}|\partial_y\Tilde{\epsilon}|^2\mu+|\Tilde{\epsilon}|^2\mu-\int \left|\frac{\lambda}{\alpha}y+1\right|^{-\sigma}\xi_aQ^2(3\Tilde{\epsilon}_1^2+\Tilde{\epsilon}_2^2)\mu\notag\\
    \gtrsim&\int_{|y|\leq\frac{1}{\sqrt{a}}}|\partial_y\Tilde{\epsilon}|^2\mu+|\Tilde{\epsilon}|^2\mu+\mathcal{O}(\delta^{C_k}\|\epsilon\|_{H^1_\mu}^2)\notag\\
    \gtrsim&\int_{|y|\leq\frac{1}{\sqrt{b}}}(|\partial_y{\epsilon}|^2+|{\epsilon}|^2)\mu+\mathcal{O}(\delta^{C_k}\|\epsilon\|_{H^1_\mu}^2).
\end{align}
Injecting \eqref{coer:4-1} into \eqref{coer:3-1}, we can obtain the lower bound of \eqref{Bound:J-1}.

For the upper bound, by Sobolev embedding and the smallness of $\epsilon$, we have the rough upper bound

\[|J(\Tilde{u})|\lesssim \|\nabla\Tilde{u}\|_{L^2}^2+\frac{1}{\lambda^2}\|\Tilde{u}\|_{L^2}^2.\]
Now we complete the proof of this Lemma.
\end{proof}

Next lemma we will give the upper bound of $K(\Tilde{u})$ (see \eqref{def:Ku-1}).
\begin{lemma}
We have the bound
\begin{align}\label{esti:K-1}
    |K(\Tilde{u})|\lesssim\frac{a}{\lambda^4}\|\epsilon\|_{H^1_\mu}^2.
\end{align}
\end{lemma}
\begin{proof}
By using \eqref{def:Qtilde-1}, we get
\begin{align*}
    \partial_t\Tilde{Q}=&i\gamma_t\Tilde{Q}-\frac{\lambda_t}{\lambda}\Tilde{Q}-\frac{r-\alpha(t)}{\lambda}\frac{\lambda_t}{\lambda}\frac{1}{\alpha^{\frac{\sigma}{2}}}\frac{1}{\lambda}\mathbf{Q}^{\prime}_{\mathcal{P}}\left(\frac{r-\alpha(t)}{\lambda}\right)e^{i\gamma}-\frac{\alpha_t}{\lambda}\frac{1}{\alpha^{\frac{\sigma}{2}}}\frac{1}{\lambda}\mathbf{Q}^{\prime}_{\mathcal{P}}\left(\frac{r-\alpha(t)}{\lambda}\right)e^{i\gamma}\\
    &+a_t\frac{1}{\alpha^{\frac{\sigma}{2}}}\frac{1}{\lambda}\partial_a\mathbf{Q}_{\mathcal{P}}\left(\frac{r-\alpha(t)}{\lambda}\right)e^{i\gamma}+\Tilde{\beta}_t\frac{1}{\alpha^{\frac{\sigma}{2}}}\frac{1}{\lambda}\partial_{\Tilde{\beta}}\mathbf{Q}_{\mathcal{P}}\left(\frac{r-\alpha(t)}{\lambda}\right)e^{i\gamma}-\frac{\sigma}{2}\frac{\alpha_t}{\alpha}\Tilde{Q}\\
    =&\left(\frac{i(1+\Tilde{\beta}^2)}{\lambda^2}+\frac{a}{\lambda^2}\right)\Tilde{Q}+\frac{a}{\lambda}\frac{r-\alpha(t)}{\lambda}\partial_r\Tilde{Q}+\frac{2\beta}{\lambda}\partial_r\Tilde{Q}-\frac{\sigma}{2}\frac{ a_{\infty}a}{\lambda^2}\Tilde{Q}\\
    &+\frac{1}{\lambda^3}\mathcal{O}\left(\left[a^2+\mathbf{Mod}+\left|a_s+\frac{1}{2}a^2-\frac{a}{\Tilde{\beta}}\mathbf{A}_2-a\mathbf{A}_1\right|\right]\xi_a|y|^{c_l}e^{-|y|}\right)\\
    =&\frac{i(1+\beta^2)}{\lambda^2}\Tilde{Q}+\frac{2\beta}{\lambda}\partial_r\Tilde{Q}+\frac{1}{\lambda^3}\mathcal{O}(a\xi_a|y|^{c_l}e^{-|y|}),
\end{align*}
where we used the modulation estimate \eqref{mod:esti:1-1} and the decay estimate \eqref{decay:P-1}. Then we have
\begin{align*}
&-\Re\left(\partial_t\Tilde{Q},\overline{\frac{1}{|x|^{\sigma}}\left(f(\Tilde{Q}+\Tilde{u})-f(\Tilde{Q}-f^{\prime}(\Tilde{Q})\cdot\Tilde{u})\right)}\right)\\
    =&-\Re\left(\partial_t\Tilde{Q},\overline{\frac{1}{|x|^{\sigma}}\left(\Bar{\Tilde{Q}}\Tilde{u}^2+2\Tilde{Q}|\Tilde{u}|^2+|\Tilde{u}|^2\Tilde{u}\right)}\right)\\
    =&-\frac{(1+\beta^2)}{\lambda^2}\Im\int \frac{1}{|x|^{\sigma}}\left(\Bar{\Tilde{Q}}\Tilde{u}^2+2\Tilde{Q}|\Tilde{u}|^2+|\Tilde{u}|^2\Tilde{u}\right)\Bar{\Tilde{Q}}\\
    &-\frac{2\beta}{\lambda}\Re\int \frac{1}{|x|^{\sigma}}\left(\Bar{\Tilde{Q}}\Tilde{u}^2+2\Tilde{Q}|\Tilde{u}|^2+|\Tilde{u}|^2\Tilde{u}\right)\partial_r\Bar{\Tilde{Q}}\\
    &+\frac{1}{\lambda^3}\mathcal{O}\left(\int \frac{1}{|x|^{\sigma}}a\xi_a|y|^{c_l}e^{-|y|}\left|\left(\Bar{\Tilde{Q}}\Tilde{u}^2+2\Tilde{Q}|\Tilde{u}|^2+|\Tilde{u}|^2\Tilde{u}\right)\right|\right)\\
    \lesssim&-\frac{(1+\beta^2)}{\lambda^2}\Im\int \frac{1}{|x|^{\sigma}}\left(\Bar{\Tilde{Q}}\Tilde{u}^2+2\Tilde{Q}|\Tilde{u}|^2+|\Tilde{u}|^2\Tilde{u}\right)\Bar{\Tilde{Q}}\\
    &-\frac{2\beta}{\lambda}\Re\int \frac{1}{|x|^{\sigma}}\left(\Bar{\Tilde{Q}}\Tilde{u}^2+2\Tilde{Q}|\Tilde{u}|^2+|\Tilde{u}|^2\Tilde{u}\right)\partial_r\Bar{\Tilde{Q}}+\frac{a}{\lambda^4}\left(1+\|\epsilon\|_{L^\infty(y\geq-\frac{\delta}{a})}\right)\int|\epsilon|^2\mu\\
    \lesssim&-\frac{(1+\beta^2)}{\lambda^2}\Im\int \frac{1}{|x|^{\sigma}}\left(\Bar{\Tilde{Q}}\Tilde{u}^2+2\Tilde{Q}|\Tilde{u}|^2+|\Tilde{u}|^2\Tilde{u}\right)\Bar{\Tilde{Q}}\\
    &-\frac{2\beta}{\lambda}\Re\int \frac{1}{|x|^{\sigma}}\left(\Bar{\Tilde{Q}}\Tilde{u}^2+2\Tilde{Q}|\Tilde{u}|^2+|\Tilde{u}|^2\Tilde{u}\right)\partial_r\Bar{\Tilde{Q}}+\frac{a}{\lambda^4}\|\epsilon\|_{H^1_\mu}^2,
\end{align*}
where we used the priori estimate \eqref{priori:1-1}. Injecting this into \eqref{def:Ku-1}, we get
\begin{align}\label{esti:K2-1}
    K(\Tilde{u})=&-\frac{1+\beta^2}{\lambda^2}\Im\int \frac{1}{|x|^{\sigma}}\left(\Bar{\Tilde{u}}^2\Tilde{Q}^2+|\Tilde{u}|^2(2\Tilde{Q}\Bar{\Tilde{u}}+\Bar{\Tilde{Q}}\Tilde{u})\right)\notag\\
    &-\frac{2\beta}{\lambda}\Re\int \phi\left(\frac{r}{\alpha(t)}-1\right)\frac{1}{|x|^{\sigma}}(f(\Tilde{Q}+\Tilde{u})-f(\Tilde{Q}))\overline{\partial_r\Tilde{u}}\notag\\
    &-\frac{(1+\beta^2)}{\lambda^2}\Im\int \frac{1}{|x|^{\sigma}}\left(\Bar{\Tilde{Q}}\Tilde{u}^2+2\Tilde{Q}|\Tilde{u}|^2+|\Tilde{u}|^2\Tilde{u}\right)\Bar{\Tilde{Q}}\notag\\
    &-\frac{2\beta}{\lambda}\Re\int \frac{1}{|x|^{\sigma}}\left(\Bar{\Tilde{Q}}\Tilde{u}^2+2\Tilde{Q}|\Tilde{u}|^2+|\Tilde{u}|^2\Tilde{u}\right)\partial_r\Bar{\Tilde{Q}}+\mathcal{O}\left(\frac{a}{\lambda^4}\|\epsilon\|_{H^1_\mu}\right)\notag\\
    =&-\frac{2\beta}{\lambda}\Re\int \phi\left(\frac{r}{\alpha(t)}-1\right)\frac{1}{|x|^{\sigma}}(f(\Tilde{Q}+\Tilde{u})-f(\Tilde{Q}))\overline{\partial_r\Tilde{u}}\notag\\
    &-\frac{2\beta}{\lambda}\Re\int \frac{1}{|x|^{\sigma}}\left(\Bar{\Tilde{Q}}\Tilde{u}^2+2\Tilde{Q}|\Tilde{u}|^2+|\Tilde{u}|^2\Tilde{u}\right)\partial_r\Bar{\Tilde{Q}}+\mathcal{O}\left(\frac{a}{\lambda^4}\|\epsilon\|_{H^1_\mu}^2\right)\notag\\
    =&I_1+\mathcal{O}\left(\frac{a}{\lambda^4}\|\epsilon\|_{H^1_\mu}^2\right).
\end{align}
Next, we need to estimate the remain terms in \eqref{esti:K2-1}. To estimate this, we need to introduce the cutoff function $\rho$. Let $\rho$ be a smooth compactly supported cutoff function which is $1$ in the neighborhood of the support of $\phi$, and $0$ in the neighborhood of $z=-1$. Now we compute
\begin{align*}
   I_2=&-\frac{2\beta}{\lambda}\Re\int \rho\left(\frac{r}{\alpha(t)}-1\right)\frac{1}{|x|^{\sigma}}(f(\Tilde{Q}+\Tilde{u})-f(\Tilde{Q}))\overline{\partial_r\Tilde{u}}\notag\\
    &-\frac{2\beta}{\lambda}\Re\int\rho\left(\frac{r}{\alpha(t)}-1\right) \frac{1}{|x|^{\sigma}}\left(\Bar{\Tilde{Q}}\Tilde{u}^2+2\Tilde{Q}|\Tilde{u}|^2+|\Tilde{u}|^2\Tilde{u}\right)\partial_r\Bar{\Tilde{Q}}\\
    =&-\frac{2\beta}{\lambda}\Re\int \rho\left(\frac{r}{\alpha(t)}-1\right)\frac{1}{|x|^{\sigma}}f(\Tilde{Q}+\Tilde{u})(\overline{\partial_r\Tilde{Q}+\partial_r\Tilde{u}})+\frac{2\beta}{\lambda}\Re\int\rho\left(\frac{r}{\alpha(t)}-1\right)a(r)f(\Tilde{Q})\overline{\partial_r\Tilde{Q}}\\
    &+\frac{2\beta}{\lambda}\Re\int\rho\left(\frac{r}{\alpha(t)}-1\right)\frac{1}{|x|^{\sigma}}\left(f(\Tilde{Q})\overline{\partial_r\Tilde{u}}+f^\prime(\Tilde{Q})\cdot\Tilde{u}\overline{\partial_r\Tilde{Q}}\right)\\
    =&-\frac{2\beta}{\lambda}\Re\int\rho\left(\frac{r}{\alpha(t)}-1\right)\frac{1}{|x|^{\sigma}}\partial_r\left(F(u)-F(\Tilde{Q})- f(\Tilde{Q})\Bar{\Tilde{u}}\right).
    \end{align*}
Integrating by parts in $r$ and using the properties of $\rho$, we have
\begin{align}\label{estim:I3-1}
    I_2=&\frac{2\beta}{\lambda}\Re\int \frac{1}{\alpha(t)}\rho^\prime\left(\frac{r}{\alpha(t)}-1\right)\frac{1}{|x|^{\sigma}}\left(F(u)-F(\Tilde{Q})- f(\Tilde{Q})\Bar{\Tilde{u}}\right)\notag\\
    &-\frac{2\beta}{\lambda}\Re\int \rho\left(\frac{r}{\alpha(t)}-1\right)\frac{\sigma}{|x|^{\sigma+1}}\left(F(u)-F(\Tilde{Q})- f(\Tilde{Q})\Bar{\Tilde{u}}\right)\notag\\
    =&\frac{2\beta}{\lambda}\Re\int \frac{1}{\alpha(t)}\rho^\prime\left(\frac{r}{\alpha(t)}-1\right)\frac{1}{|x|^{\sigma}}\left(\frac{1}{2}|\Tilde{Q}|^2|\Tilde{u}|^2+(\Re(\Tilde{Q}\bar{\Tilde{u}}))^2+\frac{1}{4}|\Tilde{u}|^4+|\epsilon|^2\Re(\Tilde{Q}\Bar{\Tilde{u}})\right)\notag\\
    &-\frac{2\beta}{\lambda}\Re\int \rho\left(\frac{r}{\alpha(t)}-1\right)\frac{\sigma}{|x|^{\sigma+1}}\left(\frac{1}{2}|\Tilde{Q}|^2|\Tilde{u}|^2+(\Re(\Tilde{Q}\bar{\Tilde{u}}))^2+\frac{1}{4}|\Tilde{u}|^4+|\epsilon|^2\Re(\Tilde{Q}\Bar{\Tilde{u}})\right)\notag\\
    =&\mathcal{O}\left(\frac{a}{\lambda^4}\|\epsilon\|_{H^1_\mu}^2\right),
\end{align}
where we used the properties of $\rho$,
\begin{align*}
    \frac{1}{r}\sim\frac{1}{\alpha(t)}~~\text{on the support of}~~\rho\left(\frac{\cdot}{\alpha(t)}-1\right).
\end{align*}
On the other hand, since $\rho=1$ on the support of $\phi$, we have
\begin{align*}
    &\frac{2\title{\beta}}{\lambda}\Re\int\phi\left(\frac{r}{\alpha(t)}-1\right)\frac{1}{|x|^{\sigma}}(f(\Tilde{Q}+\Tilde{u})-f(\Tilde{Q}))\partial_r\Bar{\Tilde{u}}\\
    =&\frac{2\title{\beta}}{\lambda}\Re\int(\rho\phi)\left(\frac{r}{\alpha(t)}-1\right)\frac{1}{|x|^{\sigma}}(f(\Tilde{Q}+\Tilde{u})-f(\Tilde{Q}))\partial_r\Bar{\Tilde{u}}.
\end{align*}
Thus
\begin{align*}
    I_3=&-\frac{2\title{\beta}}{\lambda}\Re\int\phi\left(\frac{r}{\alpha(t)}-1\right)\frac{1}{|x|^{\sigma}}(f(\Tilde{Q}+\Tilde{u})-f(\Tilde{Q}))\partial_r\Bar{\Tilde{u}}\notag\\
    &+\frac{2\title{\beta}}{\lambda}\Re\int\rho\left(\frac{r}{\alpha(t)}-1\right)\frac{1}{|x|^{\sigma}}(f(\Tilde{Q}+\Tilde{u})-f(\Tilde{Q}))\partial_r\Bar{\Tilde{u}}\notag\\
    =&-\frac{2\title{\beta}}{\lambda}\Re\int\rho\left(\frac{r}{\alpha(t)}-1\right)\left[\phi\left(\frac{r}{\alpha(t)}-1\right)-1\right]\frac{1}{|x|^{\sigma}}(f(\Tilde{Q}+\Tilde{u})-f(\Tilde{Q}))\partial_r\Bar{\Tilde{u}}\notag\\
    =&-\frac{2\title{\beta}}{\lambda}\Re\int\rho\left(\frac{r}{\alpha(t)}-1\right)\left[\phi\left(\frac{r}{\alpha(t)}-1\right)-1\right]\frac{1}{|x|^{\sigma}}\partial_r\left[F(\Tilde{Q}+\Tilde{u})-F(\Tilde{Q})-f(\Tilde{Q})\Bar{\Tilde{u}}\right]\notag\\
    &+\frac{2\title{\beta}}{\lambda}\Re\int\rho\left(\frac{r}{\alpha(t)}-1\right)\left[\phi\left(\frac{r}{\alpha(t)}-1\right)-1\right]\frac{1}{|x|^{\sigma}}\left(2|\Tilde{u}|^2\Tilde{Q}+\Tilde{u}^2\Bar{\Tilde{Q}}+|\Tilde{u}|^2\Tilde{u}\right)\partial_r\Bar{\Tilde{Q}}\notag\\
    =&\frac{2\beta}{\lambda}\Re\int\left(\frac{ a_{\infty}a}{2\beta\lambda}\partial_z(\rho(z)(\phi(z)-1))\frac{1}{|x|^{\sigma}}+\rho(z)(\phi(z)-1)\partial_r\frac{1}{|x|^{\sigma}}\right)\left[F(\Tilde{Q}+\Tilde{u})-F(\Tilde{Q})-f(\Tilde{Q})\Bar{\Tilde{u}}\right]\notag\\
    &+\frac{2\title{\beta}}{\lambda}\Re\int\rho\left(\frac{r}{\alpha(t)}-1\right)\left[\phi\left(\frac{r}{\alpha(t)}-1\right)-1\right]\frac{1}{|x|^{\sigma}}\left(2|\Tilde{u}|^2\Tilde{Q}+\Tilde{u}^2\Bar{\Tilde{Q}}+|\Tilde{u}|^2\Tilde{u}\right)\partial_r\Bar{\Tilde{Q}},
\end{align*}
where $z=\frac{r}{\alpha(t)}-1=\frac{ a_{\infty}a}{2\beta}y$.
Since $\phi(0)=1$, from \eqref{def:a-1}, we have
\begin{align*}
    |\phi(z)-1|\lesssim|z|\lesssim a|y|,~~~|\partial_z(\phi(z)-1)|\lesssim1.
\end{align*}
Then, we deduce
\begin{align}\label{estim:I4-1}
    I_3\lesssim&\frac{a}{\lambda^4}\int \left(|\epsilon|^4+a|y|^{c_l}\xi_ae^{-2|y|}|\epsilon|^2\right)\mu+\frac{a}{\lambda^4}\int |y|\left(|\epsilon|^2|y|^{c_l}\xi_ae^{-2|y|}+\xi_a|\epsilon|^3e^{-c|y|}\right)\mu\notag\\
    \lesssim&\frac{a}{\lambda^4}\|\epsilon\|_{H^1_{\mu}}^2,
\end{align}
where we used the Sobolev inequality \eqref{GN-1} and the decay estimate \eqref{decay:P-1}. From \eqref{estim:I3-1} and \eqref{estim:I4-1}, we get
\begin{align}\label{relation:I3:4-1}
    I_2-I_3=\mathcal{O}\left(\frac{a}{\lambda^4}\|\epsilon\|_{H^1_{\mu}}^2\right).
\end{align}
Notice that the function $1-\rho$ is supported by construction in $y\leq-\frac{1}{a}$  where $\Tilde{Q}$ vanishes, hence from \eqref{relation:I3:4-1}, we get
\begin{align}\label{esti:K3-1}
    I_1=\mathcal{O}\left(\frac{a}{\lambda^4}\|\epsilon\|_{H^1_{\mu}}^2\right).
\end{align}
Injecting \eqref{esti:K3-1} into \eqref{esti:K2-1}, we get
\[K(\Tilde{u})=\mathcal{O}\left(\frac{a}{\lambda^4}\|\epsilon\|_{H^1_{\mu}}^2\right).\]
This means that \eqref{esti:K-1} holds and we complete the proof of this lemma.
\end{proof}

\section{Backwards propagation of smallness}

In this section, we now apply the energy estimate of the previous section in order to establish a bootstrap argument that will be needed in the construction of  the ring blowup solution.

From now on, we choose the integer $l$ appearing in Lemma \ref{lemma:app:pro-1}  such that $l>5$.
Given $t_0<t_1<0$, let $u(t)$ be the solution to \eqref{equ:3-1}  with initial data at $t=t_1$ given explicitly by
\begin{align}\label{B:1-1}
    u(t_1,r)=\alpha_1^{\frac{\sigma}{2}}(t_1)\frac{1}{\lambda_1(t_1)}\mathbf{Q}_{(a_1(t_1),\Tilde{\beta}_1(t_1))}\left(\frac{r-\alpha_1(t_1)}{\lambda_1(t_1)}\right)e^{i\gamma_1(t_1)}.
\end{align}
Our aim is to derive bounds on $u$ backward on a time interval independent of $t_1$ as $t_1\to0$. From \eqref{B:1-1}, we have the well-prepared data initialization
\begin{align*}
    \epsilon(t_1)=0,~~(\lambda,a,\Tilde{\beta},\alpha,\gamma)(t_1)=(\lambda_1,a_1,\Tilde{\beta}_1,\alpha_1,\gamma_1)(t_1),
\end{align*}
and we may thus consider a backward time $t<t_1$ such that the following bootstrap assumption
\begin{align*}
    \|\epsilon\|_{H^1_{\mu}}<\min\{a,\lambda\}\delta,~~0<a<\delta,~~|\Tilde{\beta}|\leq a,\\
    \frac{b_{\infty}}{2}\leq\frac{\alpha(t)}{\lambda^{1-\frac{a_\infty}{1-\sigma a_\infty}}}\leq 2b_{\infty}.
\end{align*}
Now we claim the following estimates hold.
\begin{lemma}\label{lemma:boot-1}
Let $t_0$ be defined in Lemma \ref{lemma:initial-1}. For any $t\in[t_0,t_1)$, $t_1<0$, such that the priori estimates \eqref{priori:1-1} and \eqref{priori:2-1} are satisfied on the interval $[t_0,t_1)$, there holds
\begin{align*} &\|\epsilon\|_{H^1_\mu}\lesssim\min\{|t|^{\frac{1-(\sigma+1) a_{\infty}}{1-(\sigma-1) a_{\infty}}},\lambda\}|t|^{{\frac{1-(\sigma+1) a_{\infty}}{1-(\sigma-1) a_{\infty}}}},\notag\\
&a(t)=\frac{1}{(1-(2\sigma+1) a_{\infty})}B_1^{\frac{1-(\sigma+1) a_{\infty}}{1-(\sigma-1) a_{\infty}}}B_2^{-2\frac{1-\sigma a_{\infty}}{1-(\sigma-1) a_{\infty}}}|t|^{\frac{1-(\sigma+1) a_{\infty}}{1-(\sigma-1) a_{\infty}}}\left(1+\mathcal{O}\left(\log|t||t|^{\frac{1-(\sigma+1) a_{\infty}}{1-(\sigma-1) a_{\infty}}}\right)\right),\\ &|\Tilde{\beta}|=\mathcal{O}\left(|t|^{\frac{1-(\sigma+1)a_{\infty}}{1-(\sigma-1)a_{\infty}}}\right),\\
    &\lambda(t)=B_2^{-\frac{1}{1+\frac{ a_{\infty}}{1-\sigma a_{\infty}}}} B_1^{\frac{1}{1+\frac{a_{\infty}}{1-\sigma a_{\infty}}}}|t|^{\frac{1}{1+\frac{a_{\infty}}{1-\sigma a_{\infty}}}}\left(1+\mathcal{O}\left(\log|t||t|^{\frac{1-(\sigma+1) a_{\infty}}{1-(\sigma-1) a_{\infty}}}\right)\right),\\
    &\alpha(t)=b\lambda^{\frac{ a_{\infty}}{1-\sigma a_{\infty}}}=b_{\infty}B_2^{\frac{\frac{a_{\infty}}{1-\sigma a_{\infty}}}{1+\frac{a_{\infty}}{1-\sigma a_{\infty}}}}B_1^{\frac{\frac{a_{\infty}}{1-\sigma a_{\infty}}}{1+\frac{a_{\infty}}{1-\sigma a_{\infty}}}}|t|^{\frac{\frac{a_{\infty}}{1-\sigma a_{\infty}}}{1+\frac{a_{\infty}}{1-\sigma a_{\infty}}}}\left(1+\mathcal{O}\left(\log|t||t|^{\frac{1-(\sigma+1) a_{\infty}}{1-(\sigma-1) a_{\infty}}}\right)\right),
\end{align*}
where $B_1$, $B_2$ and $b_{\infty}$ are defined in Lemma \ref{lemma:initial-1}.
\end{lemma}
\begin{proof}
We  divide the proof into the two steps.

{\bf Step 1:} Control of $\epsilon$. By \eqref{energy:virial-1} and \eqref{def:Ku-1}, we have
\begin{align*}
    \frac{d}{dt}\frac{J(\Tilde{u})}{\lambda^m}=&\frac{1}{\lambda^m}J(\Tilde{u})+m\frac{a}{\lambda^{2+m}}J(\Tilde{u})-m\frac{\mathbf{A}_1}{\lambda^{2+m}}J(\Tilde{u})\\
    &-\left(\frac{\lambda_s}{\lambda}+a-\mathbf{A}_1\right)\frac{J(\Tilde{u})}{\lambda^{2+m}}+\mathcal{O}\left(\frac{a}{\lambda^{4+m}}\|\epsilon\|_{H^1_\mu}^2+\frac{a^{l}}{\lambda^{4+m}}\|\epsilon\|_{H^1_\mu}\right).
\end{align*}
From \eqref{mod:esti:1-1}, \eqref{priori:1-1} and \eqref{Bound:J-1}, we get
\begin{align*}
    \left|-m\frac{\mathbf{A}_1}{\lambda^{2+m}}J(\Tilde{u})\right|-\left|\left(\frac{\lambda_s}{\lambda}-a+\mathbf{A}_1\right)\frac{J(\Tilde{u})}{\lambda^{2+m}}\right|&\lesssim\frac{a}{\lambda^{4+m}}\left(a^2+a\|\epsilon\|_{H^1_\mu}+a^l\right)\|\epsilon\|_{H^1_\mu}^2\\
    &\lesssim\frac{a\delta^C}{\lambda^{4+m}}\|\epsilon\|_{H^1_\mu}^2.
\end{align*}
Again, using \eqref{Bound:J-1}, there exists $C>0$ such that
\begin{align*}
    \frac{d}{dt}\frac{J(\Tilde{u})}{\lambda^m}\geq (c_0m-C)\frac{a}{\lambda^{4+m}}\|\epsilon\|_{H^1_{\mu}}^2-C\frac{a^{2l-1}}{\lambda^{4+m}}.
\end{align*}
Then, if $c_0m-C>0$,
\begin{align}\label{Boot:1-1}
    \frac{d}{dt}\frac{J(\Tilde{u})}{\lambda^m}\gtrsim -\frac{a^{2k-1}}{\lambda^{4+m}}\gtrsim-a\lambda^{2(l-1)\left(1-\frac{ a_{\infty}}{1-\sigma a_{\infty}}\right)-4-m},
\end{align}
where in the last step we used $a\sim\lambda^{1-\frac{ a_{\infty}}{1-\sigma a_{\infty}}}$, this can be obtained by the definition of \eqref{def:a-1}, priori estimates \eqref{priori:1-1} and \eqref{priori:2-1}, and Lemma \ref{lemma:initial-1}.
Integrating from $t_0$ to $t_1$, and using $J(\Tilde{u}(t_1))=0$, we obtain
\begin{align*}
    J(\Tilde{u})\lesssim\lambda^m\int_{t_0}^{t_1}a(\tau)\lambda(\tau)^{2(l-1)\left(1-\frac{ a_{\infty}}{1-\sigma a_{\infty}}\right)-4-m}d\tau.
\end{align*}
From \eqref{mod:esti:1-1}, \eqref{priori:1-1} and \eqref{priori:2-1}, we deduce
\[\left|\frac{\lambda_s}{\lambda}+a\right|\lesssim|\mathbf{A}_1|+a\|\epsilon\|_{H^1_\mu}+a^l\lesssim\delta^{C_l}a.\]
This means that $0<a\lesssim-\lambda\lambda_t$. Therefore, if  choose $l$ is sufficient large, we have
\begin{align*}
    J(\Tilde{u})\lesssim\lambda^2(t).
\end{align*}
By \eqref{Bound:J-1}, we have
\begin{align}\label{Boot:2-1}
    \|\nabla\Tilde{u}\|_{L^2}+\frac{\|\Tilde{u}\|_{L^2}^2}{\lambda^2}\lesssim\lambda^2(t).
\end{align}
This is equivalent to
\begin{align*}
    \|\epsilon\|_{H^1_{\mu}}^2\lesssim\lambda^4(t).
\end{align*}

{\bf Step 2:} Control of the modulation parameters.  From \eqref{Boot:1-1} and \eqref{Boot:2-1}, we can obtain
\begin{align*}
    \|\epsilon\|_{H^1_{\mu}}\lesssim\lambda^{2+(l-1)\left(1-\frac{ a_{\infty}}{1-\sigma a_{\infty}}-\frac{2}{2(l-1)
    }\right)}.
\end{align*}
Together with \eqref{mod:esti:1-1} and $a\sim\lambda^{1-\frac{ a_{\infty}}{1-\sigma a_{\infty}}}$, we have
\begin{align*}
    \mathbf{Mod}(t)\lesssim a\|\epsilon\|_{H^1_\mu}+a^l\lesssim a^l.
\end{align*}
 Let $l>\frac{2}{1-\frac{a_{\infty}}{1-\sigma a_{\infty}}}+1$ and $t_1$ be as defined in Lemma \ref{lemma:initial-1} and $t_0\leq t<t_1<0$.  Let $(\lambda_0,a_0,\Tilde{\beta}_0, \alpha_{0},\gamma_0)$ be the solution to the system \eqref{dynamical:sys-1}. Let  $(\lambda,a,\Tilde{\beta},\alpha,\gamma)(t_1)$ be the initial data as
\begin{align*}
    (\lambda,a,\Tilde{\beta},\alpha,\gamma)(t_1)=(\lambda_0,a_0,\Tilde{\beta}_0, \alpha_0,\gamma_0)(t_1)
\end{align*}
and be the solution of the following perturbed system of modulation equations on $[t_0,t_1)$:
\begin{align*}
    \begin{cases}
    \frac{\lambda_s}{\lambda}+a-\mathbf{A}_1=\mathcal{O}(a^l),\\
    \frac{\alpha_s}{\lambda}+2\beta=\mathcal{O}(a^l),\\
    \Tilde{\beta}_s-\mathbf{A}_2=\mathcal{O}(a^l),\\
    a=\frac{2\beta}{ a_{\infty}}\frac{\lambda}{\alpha},~~\beta=\beta_{\infty}+\Tilde{\beta},\\
    \gamma_s=1+\beta^2+\mathcal{O}(a^l).
    \end{cases}
\end{align*}
Then we claim the following bounds hold on the time interval $[t_0,t_1)$:
\begin{align*}
&a(t)=\frac{1}{(1-(2\sigma+1) a_{\infty})}B_1^{\frac{1-(\sigma+1) a_{\infty}}{1-(\sigma-1) a_{\infty}}}B_2^{-2\frac{1-\sigma a_{\infty}}{1-(\sigma-1) a_{\infty}}}|t|^{\frac{1-(\sigma+1) a_{\infty}}{1-(\sigma-1) a_{\infty}}}\left(1+\mathcal{O}\left(\log|t||t|^{\frac{1-(\sigma+1) a_{\infty}}{1-(\sigma-1) a_{\infty}}}\right)\right),\\
&|\Tilde{\beta}|=\mathcal{O}\left(|t|^{\frac{1-(\sigma+1)a_{\infty}}{1-(\sigma-1)a_{\infty}}}\right)\\
    &\lambda(t)=B_2^{-\frac{1}{1+\frac{ a_{\infty}}{1-\sigma a_{\infty}}}} B_1^{\frac{1}{1+\frac{a_{\infty}}{1-\sigma a_{\infty}}}}|t|^{\frac{1}{1+\frac{a_{\infty}}{1-\sigma a_{\infty}}}}\left(1+\mathcal{O}\left(\log|t||t|^{\frac{1-(\sigma+1) a_{\infty}}{1-(\sigma-1) a_{\infty}}}\right)\right),\\
    &\alpha(t)=b\lambda^{\frac{ a_{\infty}}{1-\sigma a_{\infty}}}=b_{\infty}B_2^{\frac{\frac{a_{\infty}}{1-\sigma a_{\infty}}}{1+\frac{a_{\infty}}{1-\sigma a_{\infty}}}}B_1^{\frac{\frac{a_{\infty}}{1-\sigma a_{\infty}}}{1+\frac{a_{\infty}}{1-\sigma a_{\infty}}}}|t|^{\frac{\frac{a_{\infty}}{1-\sigma a_{\infty}}}{1+\frac{a_{\infty}}{1-\sigma a_{\infty}}}}\left(1+\mathcal{O}\left(\log|t||t|^{\frac{1-(\sigma+1) a_{\infty}}{1-(\sigma-1) a_{\infty}}}\right)\right),\\
    &\gamma(t)=(1+\beta_{\infty}^2)B_1^{-\frac{1-\frac{a_{\infty}}{1-\sigma a_{\infty}}}{1+\frac{a_{\infty}}{1-\sigma a_{\infty}}}}B_2^{2\frac{1}{1+\frac{a_{\infty}}{1-\sigma a_{\infty}}}}|t|^{-\frac{1-\frac{a_{\infty}}{1-\sigma a_{\infty}}}{1+\frac{a_{\infty}}{1-\sigma a_{\infty}}}}+\mathcal{O}(|\log t|),
\end{align*}
where
\begin{align*}
    B_1=\frac{1-(\sigma-1) a_{\infty}}{1-(\sigma+1) a_{\infty}},~~B_2=\frac{ a_{\infty}}{2(1-3 a_{\infty})\beta_{\infty}}b_{\infty}
\end{align*}
where $b_\infty$ is defined in Lemma \ref{lemma:initial-1}. By the similar argument as  \cite{MRS2014Duke}, we can obtain the above estimates.
This complete the proof of this lemma.
\end{proof}
%{\color{red}\noindent\rule[0.8\baselineskip]{\textwidth}{1pt}}

\section{Existence of the ring blowup solution}
In this section, we aim to prove our main result Theorem \ref{thm-1}.

\begin{proof}[\bf Proof of Theorem \ref{thm-1}.]
Let $\{t_n\}_{n\geq0}$ be an increasing sequence of times $t_n<0$ such that $t_n\to0^-.$ Let $u_n$ be the solution to \eqref{equ:3-1} with the initial data at $t=t_n$,
\begin{align*}
    u_n(t_n,x)=( \alpha_0(t_n))^{\frac{\sigma}{2}}\frac{1}{\lambda_0(t_n)}\mathbf{Q}_{\mathcal{P}_0(t_n)}\left(\frac{r- \alpha_0(t_n)}{\lambda_0(t_n)}\right)e^{i\gamma_0(t_n)},
\end{align*}
with $\mathcal{P}_0(t_n)=(a_0(t_n), \Tilde{\beta}_0(t_n))$.
Let $t_1<0$ (see Lemma \ref{lemma:boot-1}) which is independent of $n$. The $L^2$ compactness of $u_n(t_0)$ is a consequence of a standard localization procedure. Indeed, Lemma \ref{lemma:boot-1} ensures the uniform bound $\|u_n\|_{H^1}\lesssim1$.  We note that the uniform bound
\[\left|\frac{d}{dt}\int\chi_R|u_n|^2\right|=2\left|\Re\int(\nabla\chi_R\cdot\nabla u_n)\Bar{u}_n\right|\lesssim\frac{1}{R},\]
with a smooth cut-off function $\chi_R(x)=\chi\left(\frac{x}{R}\right)$ where $\chi(x)\equiv0$ for $|x|\leq1$ and $\chi(x)=1$ for $|x|\geq2$. By integrating this bound from $t_0$ to $t_1$, we have
\[\lim_{R\to\infty}\int_{|x|\geq R}|u_n(t_0)|^2=0,\]
which together with the $L^2(|x|<R)$ compactness of $u_n(t_0)$, up to a subsequence, we can obtain
\[u_n(t_0)\to u(t_0)~~~\text{in}~~L^2~~~\text{as}~~n\to+\infty.\]
Let $u\in C([t_0,0),H^1)$ be the solution to \eqref{equ:3-1} with the initial data $u(t_0)$, then using the uniform estimate in $H^1$  for $u_n$ and the converge in $L^2$ of $u_n(t_0)$, we have, for $t\in[t_0,0)$,
\[u_n(t)\to u(t)~~\text{as}~~L^2.\]
Let $u_n(t)$ admits a geometrical decomposition
\[u_n(t,x)=(\alpha_n(t))^{\frac{\sigma}{2}}\frac{1}{\lambda_n(t)}[\mathbf{Q}_{\mathcal{P}_n(t)}+\epsilon_n]\left(\frac{r-\alpha_n(t)}{\lambda_n(t)}\right)e^{i\gamma_n(t)},\]
then $u$ satisfied the following form
\[u(t,x)=\alpha(t)^{\frac{\sigma}{2}}\frac{1}{\lambda(t)}[\mathbf{Q}_{\mathcal{P}(t)}+\epsilon]\left(\frac{r-\alpha(t)}{\lambda(t)}\right)e^{i\gamma(t)},
\]
with $\mathcal{P}_n\to\mathcal{P}$, $\gamma_n\to\gamma$ and $\epsilon_n\to\epsilon$ in $L^2$ as $n\to\infty$. In addition, by Lemma \ref{lemma:boot-1}, the parameters and $\epsilon$ have the following bounds,
\begin{align*}
&a(t)=\frac{1}{(1-(2\sigma+1) a_{\infty})}B_1^{\frac{1-(\sigma+1) a_{\infty}}{1-(\sigma-1) a_{\infty}}}B_2^{-2\frac{1-\sigma a_{\infty}}{1-(\sigma-1) a_{\infty}}}|t|^{\frac{1-(\sigma+1) a_{\infty}}{1-(\sigma-1) a_{\infty}}}\left(1+\mathcal{O}\left(\log|t||t|^{\frac{1-(\sigma+1) a_{\infty}}{1-(\sigma-1) a_{\infty}}}\right)\right),\\
&|\Tilde{\beta}|=\mathcal{O}\left(|t|^{\frac{1-(\sigma+1)a_{\infty}}{1-(\sigma-1)a_{\infty}}}\right)\\
    &\lambda(t)=B_2^{-\frac{1}{1+\frac{ a_{\infty}}{1-\sigma a_{\infty}}}} B_1^{\frac{1}{1+\frac{a_{\infty}}{1-\sigma a_{\infty}}}}|t|^{\frac{1}{1+\frac{a_{\infty}}{1-\sigma a_{\infty}}}}\left(1+\mathcal{O}\left(\log|t||t|^{\frac{1-(\sigma+1) a_{\infty}}{1-(\sigma-1) a_{\infty}}}\right)\right),\\
    &\alpha(t)=b\lambda^{\frac{ a_{\infty}}{1-\sigma a_{\infty}}}=b_{\infty}B_2^{\frac{\frac{a_{\infty}}{1-\sigma a_{\infty}}}{1+\frac{a_{\infty}}{1-\sigma a_{\infty}}}}B_1^{\frac{\frac{a_{\infty}}{1-\sigma a_{\infty}}}{1+\frac{a_{\infty}}{1-\sigma a_{\infty}}}}|t|^{\frac{\frac{a_{\infty}}{1-\sigma a_{\infty}}}{1+\frac{a_{\infty}}{1-\sigma a_{\infty}}}}\left(1+\mathcal{O}\left(\log|t||t|^{\frac{1-(\sigma+1) a_{\infty}}{1-(\sigma-1) a_{\infty}}}\right)\right),\\
    &\gamma(t)=(1+\beta_{\infty}^2)B_1^{-\frac{1-\frac{a_{\infty}}{1-\sigma a_{\infty}}}{1+\frac{a_{\infty}}{1-\sigma a_{\infty}}}}B_2^{2\frac{1}{1+\frac{a_{\infty}}{1-\sigma a_{\infty}}}}|t|^{-\frac{1-\frac{a_{\infty}}{1-\sigma a_{\infty}}}{1+\frac{a_{\infty}}{1-\sigma a_{\infty}}}}+\mathcal{O}(|\log t|),
\end{align*}
and
\begin{align*}
    \|\epsilon\|_{H^1_{\mu}}\lesssim|t|^{\frac{2}{1+\frac{a_{\infty}}{1-\sigma a_{\infty}}}},
\end{align*}
where
\begin{align*}
    B_1=\frac{1-(\sigma-1) a_{\infty}}{1-(\sigma+1) a_{\infty}},~~B_2=\frac{ a_{\infty}}{2(1-3 a_{\infty})\beta_{\infty}}b_{\infty},
\end{align*}
and $b_{\infty}$ given by Lemma \ref{lemma:initial-1}.
Now by the standard argument as \cite{MRS2014Duke,RS2011JAMS,GL2022CPDE,GL2021JFA}, we can obtain $u\in C([0,T),H^1)$  and $u$ blows up at time $T=0$. The parameters estimates are the consequence of the above estimates for $(a,\lambda,\alpha,\gamma,\epsilon)$.

This complete the proof of Theorem \ref{thm-1}.
\end{proof}

\appendix

\section{Appendix}\label{appendix:unique-1}
In this section, we aim to prove Lemma \ref{lemma:decom:unique-1}. This is a standard modulation lemma relies on the implicit function theorem and the mass subcritical nondegeneracy $(Q,\Lambda Q)=\frac{1}{2}\|Q\|_{L^2}^2\neq0$.
\begin{proof}[\bf Proof of Lemma \ref{lemma:decom:unique-1}.]
By the assumption, we have
\begin{align*}
    u(x)= \alpha_0^{\frac{\sigma}{2}}(t)\frac{1}{\lambda_0(t)}Q_{(a_0,\Tilde{\beta}_0)}\left(\frac{r-\alpha_0}{\lambda_0}\right)e^{i\gamma_0}+\Tilde{u}_0(x),
\end{align*}
and we wish to introduce a modified decomposition
\begin{align*}
    u(x)=\alpha_1^{\frac{\sigma}{2}}(t)\frac{1}{\lambda_1(t)}Q_{(a_1,\Tilde{\beta}_1)}\left(\frac{r-\alpha_1}{\lambda_1}\right)e^{i\gamma_1}+\Tilde{u}_1(x).
\end{align*}
Comparing the decompositions, we obtain the formula
\begin{align*}
    \Tilde{u}_1(x)= \alpha_0^{\frac{\sigma}{2}}(t)\frac{1}{\lambda_0(t)}Q_{(a_0,\Tilde{\beta}_0)}\left(\frac{r-\alpha_0}{\lambda_0}\right)e^{i\gamma_0}-\alpha_1^{\frac{\sigma}{2}}(t)\frac{1}{\lambda_1(t)}Q_{(a_1,\Tilde{\beta}_1)}\left(\frac{r-\alpha_1}{\lambda_1}\right)e^{i\gamma_1}+\Tilde{u}_0(x).
\end{align*}
We now define the functional
\begin{align*}
    F_{z,\mu,\gamma,\Tilde{\beta
    },\nu}(y)=\nu^{\frac{\sigma}{2}}\mu Q_{a_0,\Tilde{\beta}_0}(\mu y+z)e^{i\gamma+i(\beta_0+\Tilde{\beta})y}-Q_{a_1,\Tilde{\beta}_1}(y)e^{i\Tilde{\beta}_1y},
\end{align*}
with
\begin{align*}
    z=\frac{\alpha_1- \alpha_0}{\lambda_0},~~\mu=\frac{\lambda_1}{\lambda_0},~~\gamma=\gamma_1-\gamma_0,~~\Tilde{\beta}=\Tilde{\beta}_1-\Tilde{\beta}_0,~~\nu=\frac{ \alpha_0}{\alpha_1}.
\end{align*}
So that
\begin{align*}
    \Tilde{\epsilon}_1(y)= F_{z,\mu,\gamma,\Tilde{\beta
    },\nu}(y)+\alpha_1^{\frac{\sigma}{2}}\lambda_1\Tilde{u}_0(\lambda_1y+\alpha_1)e^{i\gamma_1+i\beta_1y}.
\end{align*}
We then define the scalar products, for $j=1,2$
\begin{align*}
    \rho^{(j)}=&\int_{-\infty}^{+\infty}\Re \Tilde{\epsilon}_1\xi_{a_1}T^{(j)}dy\\
    =&\int_{-\infty}^{+\infty}\Re{F_{z,\mu,\gamma,\Tilde{\beta
    },\nu}(y)}\xi_{a_1}T^{(j)}dy+\Re\int\Tilde{u}_0(r)\alpha_1^{\frac{\sigma}{2}}(\xi_{a_1}T^{(j)})\left(\frac{r-\alpha_1}{\lambda_1}\right)e^{i\gamma-1+i\beta_1\frac{r-\alpha_1}{\lambda_1}}dr,
\end{align*}
and for $j=3,4$
\begin{align*}
    \rho^{(j)}=&\int_{-\infty}^{+\infty}\Im \Tilde{\epsilon}_1\xi_{a_1}T^{(j)}dy\\
    =&\int_{-\infty}^{+\infty}\Im{F_{z,\mu,\gamma,\Tilde{\beta
    },\nu}(y)}\xi_{a_1}T^{(j)}dy+\Re\int\Tilde{u}_0(r)\alpha_1^{\frac{\sigma}{2}}(\xi_{a_1}T^{(j)})\left(\frac{r-\alpha_1}{\lambda_1}\right)e^{i\gamma-1+i\beta_1\frac{r-\alpha_1}{\lambda_1}}dr,
\end{align*}
where
\begin{align*}
    T^{(1)}=yQ,~~T^{(2)}=Q,~~,T^{(3)}=\partial_yQ,~~T^{(4)}=\Lambda Q.
\end{align*}
We now view $\rho=(\rho^{(j)})_{1\leq j\leq 4}$ as smooth functions of $(\Tilde{u}_0,z,\mu,\Tilde{\beta},\gamma)$. Observe that the bound \eqref{lower:priori-1} ensures that
\begin{align}\notag
    |\rho(\Tilde{u}_0,0,1,0,0)|\lesssim\delta.
\end{align}
Notice that
\begin{align*}
  a_1=\frac{2\beta_1}{ a_{\infty}}\frac{\lambda_1}{\alpha_1}=2(\beta_0+\Tilde{\beta})\frac{\lambda_0}{ a_{\infty} \alpha_0}\mu\frac{ \alpha_0}{\alpha_1}=\left(1+\frac{\Tilde{\beta}}{\beta_0}\right)a_0\mu\left(1+\frac{ a_{\infty}a_0}{2\beta_0}z\right)^{-1}.
 \end{align*}
Using
\begin{align*}
    Q_{\mathcal{P}}{_{|_{\mathcal{P}=0}}}=Qe^{-i\beta_{\infty}y}
\end{align*}
we obtain
\begin{align*}
    &\partial_zF|_{(z=0,\mu=1,\Tilde{\beta}=1,\gamma=0)}=Q^{\prime}-i\beta_{\infty}Q+\mathcal{O}(|a_0|+|\Tilde{\beta}_0|)e^{-c|y|},\\
    &\partial_{\mu}F|_{(z=0,\mu=1,\Tilde{\beta}=1,\gamma=0)}=\Lambda Q-i\beta_{\infty}Q+\mathcal{O}(|a_0|+|\Tilde{\beta}_0|)e^{-c|y|},\\
    &\partial_{\Tilde{\beta}}F|_{(z=0,\mu=1,\Tilde{\beta}=1,\gamma=0)}=iyQ+\mathcal{O}(|a_0|+|\Tilde{\beta}_0|)e^{-c|y|},\\
    &\partial_{\gamma}F|_{(z=0,\mu=1,\Tilde{\beta}=1,\gamma=0)}=-iQ+\mathcal{O}(|a_0|+|\Tilde{\beta}_0|)e^{-c|y|}.
\end{align*}
Then the Jacobian matrix of $\rho$ at $(z=0,\mu=1,\Tilde{\beta}=1,\gamma=0)$ is not zero from the smallness assumption \eqref{small:para-1}. The existence of the desired decomposition now follows from the implicit function theorem, and the smallness of the parameters.
\end{proof}

% \renewcommand{\proofname}{\bf Proof.}
% \section{Nondegeneracy of \texorpdfstring{$L_+$}{L+}}

%\begin{proof}[Proof of Theorem ]

% \end{proof}

\renewcommand{\proofname}{\bf Proof.}

\noindent
{\bf Acknowledgments}

Y.Li was supported by China Postdoctoral Science Foundation (No. 2021M701365) and the funding of innovating activities in Science and Technology of Hubei Province.
%%%%%%%%%%%%%%%%%%%%%%%%%%%%%%%%%%%%%%%%%%%%%%

%%%%%%%%%%%%%%%%%%%%%%%%%%%%%%%%%%%%%%%%%%%%%%

\vspace*{.5cm}

%\bibliographystyle{plain}
%\bibliography{ref}

\bigskip

%\begin{flushleft}
% Vladimir Georgiev,\\
% Dipartimento di Matematica, Universit\`{a} di Pisa, Largo B. Pontecorvo 5, 56127 Pisa, Italy\\
%  Faculty of Science and Engineering, Waseda University, 3-4-1, Okubo, Shinjuku-ku, Tokyo 169-8555, Japan\\
%  IMICBAS, Acad. Georgi Bonchev Str., Block 8, 1113 Sofia, Bulgaria
% \quad
% E-mail: georgiev@dm.unipi.it
% \end{flushleft}

\begin{flushleft}
Yuan Li,\\
School of Mathematics and Statistics, Central China Normal University, Wuhan, PR China
\quad
E-mail: yli2021@ccnu.edu.cn
\end{flushleft}

\bigskip

\medskip
\end{document}